\documentclass{amsart}

\usepackage{latexsym,amsthm,mathrsfs,amsmath,amscd,enumerate,enumitem, amscd,color,dsfont, cancel, mathtools, textcomp, float}

\usepackage[utf8x]{inputenc}

\usepackage[cspex,bbgreekl]{mathbbol}
\usepackage{amsfonts}
\usepackage{amssymb}             % AMS Math
\usepackage{graphicx}
\DeclareSymbolFontAlphabet{\mathbbl}{bbold}
\DeclareSymbolFontAlphabet{\mathbb}{AMSb}%
 \newtheorem{thm}{Theorem}[section]

\theoremstyle{definition}

 \theoremstyle{remark}

\usepackage{tikz}
\usetikzlibrary{plotmarks}

\usepackage{hyperref}
\newcommand{\supp}{\mathop{\mathrm{supp}}}

\makeatletter
\DeclareRobustCommand\widecheck[1]{{\mathpalette\@widecheck{#1}}}
 \def\@widecheck#1#2{%
   \setbox\z@\hbox{\m@th$#1#2$}%
  \setbox\tw@\hbox{\m@th$#1%
       \widehat{%
          \vrule\@width\z@\@height\ht\z@
          \vrule\@height\z@\@width\wd\z@}$}%
    \dp\tw@-\ht\z@
    \@tempdima\ht\z@ \advance\@tempdima2\ht\tw@ \divide\@tempdima\thr@@
    \setbox\tw@\hbox{%
       \raise\@tempdima\hbox{\scalebox{1}[-1]{\lower\@tempdima\box
\tw@}}}
    {\ooalign{\box\tw@ \cr \box\z@}}}
\makeatother

\usepackage{caption}
\usepackage{subcaption}
\usepackage[left=2.2cm,top=2.5cm,right=2.2cm,bottom=2.5cm]{geometry} % Document margins

\numberwithin{equation}{section}
\allowdisplaybreaks
\begin{document}

\title{$L^p$--boundedness of Stein's square functions associated to Fourier--Bessel expansions}

\author[V. Almeida]{V\'ictor Almeida}
\address{V\'ictor Almeida, Jorge J. Betancor, Lourdes Rodr\'iguez-Mesa\newline
	Departamento de An\'alisis Matem\'atico, Universidad de La Laguna,\newline
	Campus de Anchieta, Avda. Astrof\'isico S\'anchez, s/n,\newline
	38721 La Laguna (Sta. Cruz de Tenerife), Spain}
\email{valmeida@ull.es, jbetanco@ull.es, lrguez@ull.edu.es}

\author[J. J. Betancor]{Jorge J. Betancor}

\author[E. Dalmasso]{Estefan\'ia Dalmasso}
\address{Estefan\'ia Dalmasso\newline
	Instituto de Matem\'atica Aplicada del Litoral, UNL, CONICET, FIQ.\newline Colectora Ruta Nac. N° 168, Paraje El Pozo,\newline S3007ABA, Santa Fe, Argentina}
\email{edalmasso@santafe-conicet.gov.ar}

\author[L. Rodr\'iguez-Mesa]{Lourdes Rodr\'iguez-Mesa}

\thanks{This paper is partially supported by MTM2016--79436--P}

\subjclass[2010]{42B35, 42B25.}

\keywords{Fourier--Bessel expansions, Bochner--Riesz means, Stein square functions, multipliers.}

\date{\today}

\begin{abstract}
In this paper we prove $L^p$ estimates for Stein's square functions associated to Fourier--Bessel expansions. Furthermore we prove transference results for square functions from Fourier--Bessel series to Hankel transforms. Actually, these are transference results for vector--valued multipliers from discrete to continuous in the Bessel setting. As a consequence, we deduce the sharpness of the range of $p$ for the $L^p$--boundedness of Fourier--Bessel Stein's square functions from the corresponding property for Hankel--Stein square functions. Finally, we deduce $L^p$ estimates for Fourier--Bessel multipliers from that ones we have got for our Stein square functions.
\end{abstract}

\maketitle

\section{Introduction}
Our purpose in this paper is to establish $L^p$--boundedness properties for Stein's square functions associated to Fourier--Bessel expansions. These results are applied to obtain $L^p$--inequalities for certain Fourier--Bessel multipliers. 

Before stating our results, we recall other ones that motivated our study. Firstly, we consider the Bochner--Riesz mean of order $\alpha> 0$, $R_t^\alpha$, defined by
$$
\widehat{R_t^\alpha f}(\xi)=\left(1-\frac{|\xi |^2}{t^2}\right)_{+}^\alpha \widehat{f}(\xi),\quad t>0,\,\xi\in \mathbb{R}^n.
$$
Here, by $\widehat{f}$ we denote the Fourier transform and $u_+=\max\{0,u\}$, for $u\in \mathbb{R}$. The Bochner--Riesz mean is a summability method that has been studied by many authors, who were looking for convergence properties of Fourier series and integrals. Even today, there are important and difficult open problems related to these matters. One of them is the Bochner--Riesz conjecture which says that, for every $t>0$ and $1<p<\infty$, $p\neq 2$, the inequality $\|R_t^\alpha f\|_p\leq C\|f\|_p$ holds, if and only if, $\alpha>\left(n|\frac{1}{2}-\frac{1}{p}|-\frac{1}{2}\right)_+$. The necessity of this condition can be found in \cite{Fe1}. The conjecture is true when $n=2$ (see \cite{CSj}), but the problem is still open when $n\geq 3$. However, some partial results have been obtained (see \cite{Bo1, Bo2, BG, Lee}).

The square function $G_\alpha$ defined by
\[
G_\alpha (f)=\left(\int_0^\infty \left|\frac{\partial}{\partial t}R_t^\alpha f\right|^2tdt\right)^{\frac{1}{2}}
\]
was introduced by Stein \cite{St1} in order to study almost everywhere convergence for Bochner--Riesz means for Fourier series and integrals. $L^p$--boundedness properties for $G_\alpha$ imply other ones for maximal operators associated with the Bochner--Riesz means and multipliers. Carbery (\cite{Ca1}) showed how $G_\alpha$ dominates Fourier multipliers and maximal Fourier multipliers.

Note that $G_\alpha$ can be written as
\[
G_\alpha (f)=\|K_t*f\|_{L^2\left((0,\infty ), \frac{dt}{t}\right)}
\]
where 
\[
\widehat{K_t}(y)=2\alpha \frac{|y|^2}{t^2}\left(1-\frac{|y|^2}{t^2}\right)_+^{\alpha -1},\quad t>0,\, y\in \mathbb{R}^n.
 \]

Here and throughout the paper,  the expression $(1-u)_+^{\sigma}$, when $\sigma <0$, will mean
$$
(1-u)_+^{\sigma}=(1-u)^{\sigma}\chi_{(0,1)}(u),\quad u\in (0,\infty).
$$

The above expression of $G_\alpha$ says that it can be seen as a vector--valued convolution operator. By using Plancherel theorem it follows that $G_\alpha$ is bounded on $L^2(\mathbb{R}^n)$ provided that $\alpha >\frac{1}{2}$ (\cite{St1}). If $1< p\leq2$, $G_\alpha$ is bounded on $L^p(\mathbb{R}^n)$ if and only if $\alpha >n\left(\frac{1}{p}-\frac{1}{2}\right)+\frac{1}{2}$ (see \cite{IK} and \cite{Su1}). In this case the oscillation character of the kernel $K_t$ does not play any role. When $p>2$, a necessary condition for $G_\alpha$ to be bounded on $L^p(\mathbb{R}^n)$ is that $\alpha >\max\big\{n\left(\frac{1}{2}-\frac{1}{p}\right),\frac{1}{2}\big\}$ (\cite{LRS1}). It was conjectured that this condition is also sufficient. In dimension 2, the conjecture was verified by Carbery (\cite{Ca2}). When $n\geq 3$ the question remains unsolved. Some progress has been made in \cite{Ch, Ki1, LRS2, See}. In \cite{LRS2}, some weighted $L^2$-estimates for $G_\alpha$ were obtained. Carro and Domingo-Salazar \cite{CD} proved weighted $L^p$--inequalities for $G_\alpha$ showing that, when $\alpha >\frac{n+1}{2}$, it can be controlled by a finite sum of sparse operators. 

Bochner--Riesz means are radial Fourier multipliers and they can be defined by using the spectral resolution of the Laplace operator. This fact suggests defining the Bochner--Riesz means for positive self-adjoint operators. Suppose that $(X,d,\mu)$ is a metric measure space, where $d$ is the metric and $\mu$, the measure. Let $L$ be a positive self-adjoint operator on $L^2(X)$. We can write
$$
L=\int_0^\infty \lambda dE_L(\lambda),
$$
where $E_L$ represents the spectral measure associated with $L$. The Bochner--Riesz mean $R_t^{\alpha ,L}$ of order $\alpha \geq 0$ is defined by
$$
R_t^{\alpha , L}(f)=\int_0^\infty \left(1-\frac{\lambda }{t^2}\right)_+^\alpha dE_L(\lambda )f,\quad f\in L^2(X)\mbox{ and }t>0.
$$
In other words, for every $t>0$, $R_t^{\alpha ,L}=S_t^\alpha (L)$, where $S_t^\alpha (z)=(1-z/t^2)_+^\alpha $, $z>0$.

Also, the maximal operator $R_*^{\alpha ,L}$ and the Stein's square function $G_\alpha ^L$ are defined in the natural way. $L^p$--boundedness properties for $R_t^{\alpha ,L}$, for $t>0$, $R_*^{\alpha ,L}$ and $G_\alpha ^L$ have been established in \cite{CDY, CLSY, CG, SYY}.

In order to obtain these properties it is usual to assume that the operator $L$ satisfies some kind of Gaussian or Davies-Gaffney estimates and has a finite speed of propagation property.

We now define the Stein's square function associated with Fourier--Bessel expansions, which is the main object of our study. Let $\nu >-1$. By $J_\nu$ we denote the Bessel function of the first kind and order $\nu$ and by $\{s_{\nu ,j}\}_{j=1}^\infty$ we represent the sequence of positive zeros of the function $J_\nu$ where $s_{\nu ,j}<s_{\nu, j+1}$, $j\in \mathbb{N}$. We define, for every $j\in \mathbb{N}$,
$$
\phi _{\nu ,j}(x)=\frac{\sqrt{2x}J_\nu (s_{\nu ,j}x)}{|J_{\nu +1}(s_{\nu ,j})|},\quad x\in (0,1).
$$

The sequence $\{\phi _{\nu ,j}\}_{j=1}^\infty$ is a complete and orthonormal system in $L^2(0,1)$ (\cite{Ho}). If $f$ is a complex measurable function on $(0,1)$ such that $x^{\min \{\nu +\frac{1}{2}, 0\}}f\in L^1(0,1)$, we define
$$
c_{\nu ,j}(f)=\int_0^1f(x)\phi _{\nu ,j}(x)dx,\quad j\in \mathbb{N}.
$$
Note that, since $\nu >-1$, $x^{\min \{\nu +\frac{1}{2}, 0\}}f\in L^1(0,1)$, provided that $f\in L^2(0,1)$.

The Bessel operator $L_\nu$ is defined by
$$
L_\nu f=\sum_{j=1}^\infty s_{\nu ,j}^2c_{\nu ,j}(f)\phi _{\nu ,j},\quad f\in D(L_\nu),
$$
where the domain $D(L_\nu)$ of $L_\nu$ is
$$
D(L_\nu )=\left\{f\in L^2(0,1):\sum_{j=1}^\infty s_{\nu ,j}^4|c_{\nu ,j}(f)|^2<\infty \right\} .
$$
$D(L_\nu)$ is a dense subspace of $L^2(0,1)$ because the space $C_c^\infty (0,1)$ of smooth functions with compact support in $(0,1)$ is contained in $D(L_\nu)$. Moreover, for every $\phi \in C_c^\infty (0,1)$, $L_\nu \phi=\mathbb{L}_\nu \phi$, where 
$$
\mathbb{L}_\nu =-x^{-\nu -\frac{1}{2}}\frac{d}{dx}x^{2\nu +1}\frac{d}{dx}x^{-\nu -\frac{1}{2}},\mbox{  on }(0,1).
$$
$L_\nu$ is a positive self-adjoint operator. Hence, $-L_\nu$ generates a holomorphic semigroup $\{W_t^\nu \}_{t>0}$ in $L^2(0,1)$. For every $t>0$, we can write
$$
W_t^\nu (f)(x)=\int_0^1 W_t^\nu (x,y)f(y)dy,\quad f\in L^2(0,1),
$$
where 
$$
W_t^\nu (x,y)=\sum_{j=1}^\infty \phi _{\nu ,j}(x)\phi _{\nu ,j}(y)e^{-ts_{\nu ,j}^2},\quad x,y\in (0,1),\;t>0.
$$
According to \cite[Theorem 1.1]{MSZ} (see also \cite{NR}), there exists $C>0$ such that
$$
\frac{1}{C}G_t^\nu (x,y)\leq W_t^\nu (x,y)\leq CG_t^\nu (x,y),\quad x,y\in (0,1),\;t>0,
$$
where
$$
G_t^\nu (x,y)=\frac{(xy)^{\nu +\frac{1}{2}}(1+t)^{\nu +2}}{(t+xy)^{\nu +\frac{1}{2}}}\min\left\{1, \frac{(1-x)(1-y)}{t}\right\}\frac{1}{\sqrt{t}}\exp\left(-\frac{|x-y|^2}{4t}-s_{\nu ,1}^2t\right),\quad x,y\in (0,1), \;t>0.
$$
Note that if $-1<\nu <-\frac{1}{2}$ the kernel $W_t^\nu (x,y)$ has not Gaussian upper bounds, and if $\nu \geq -\frac{1}{2}$, then
$$
W_t^\nu (x,y)\leq \frac{C}{\sqrt{t}}\exp \left(-\frac{|x-y|^2}{4t}\right),\quad x,y\in (0,1),\; t>0.
$$

The Bochner--Riesz means associated with $L_\nu$ can be represented by
$$
R_t^{\alpha ,L_\nu }(f)=\sum_{j=1}^\infty \left(1-\frac{s_{\nu ,j}^2}{t^2}\right)_+^\alpha c_{\nu ,j}(f)\phi _{\nu ,j},\quad f\in L^2(0,1)\mbox{ and }t>0.
$$
Weighted $L^p$--inequalities for $L_\nu$--Bochner--Riesz means were obtained by Ciaurri and Roncal (\cite{CR1} and \cite{CR2}). $L^p$--boundedness properties for the maximal operator associated to $\{R_t^{\alpha ,L_\nu }\}_{t>0}$ were established in \cite{CR3}.

The Stein's square function associated to $L_\nu$ is given by
$$
G_\alpha ^{L_\nu}(f)=\left(\int_0^\infty \left|\frac{\partial}{\partial t}R_t^{\alpha ,L_\nu }(f)\right|^2tdt\right)^{\frac{1}{2}}.
$$

The main result of this paper is the following.

\begin{thm}\label{Theorem1.1}
There exists $C>0$ such that
\begin{equation}\label{1.1}
\frac{1}{C}\|f\|_p\leq \|G_\alpha ^{L_\nu}(f)\|_p\leq C\|f\|_p,\quad f\in L^2(0,1)\cap L^p(0,1),
\end{equation}
provided that some of the following three properties hold:
\begin{enumerate}[label=(\roman*)]
     \item $\nu \geq -\frac12$, $1<p<\infty$ and $\alpha>\max\left\{\frac1p,1-\frac1p\right\}$;
    \item $-1<\nu<-\frac12$, $\frac{2}{2\nu+3}<p\leq 2$ and $\alpha>\frac{1}{4\nu+4}\left(\frac2p+2\nu+1\right)$;
    \item $-1<\nu<-\frac12$, $2\leq p<-\frac{2}{2\nu+1}$ and $\alpha>\frac{1}{4\nu+4}\left(-\frac2p+2\nu+3\right)$.
\end{enumerate}
\end{thm}
%\newpage
The following figures illustrate the situation.

\begin{figure}[!htb]
    \centering
    \begin{minipage}{.5\textwidth}
        \centering
        \includegraphics{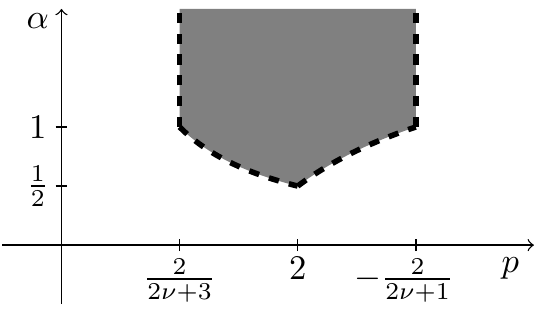}
        \caption{Case $-1<\nu<-\frac{1}{2}$}
        \label{fig:prob1_6_2}
    \end{minipage}%
    \begin{minipage}{0.5\textwidth}
        \centering
        \includegraphics{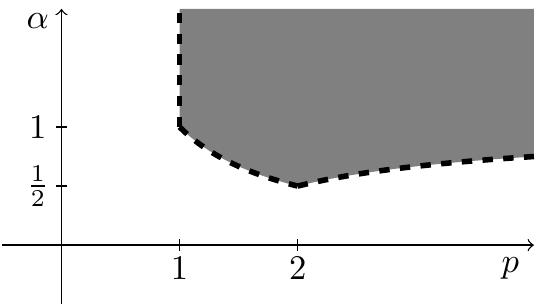}
        \caption{Case $\nu\geq - \frac{1}{2}$}
        \label{fig:prob1_6_1}
    \end{minipage}
\end{figure}

Note that if $-1<\nu <-\frac{1}{2}$ and $1<p<2$, then $\frac{1}{p}<\frac{1}{4\nu+4}\left(\frac{2}{p}+2\nu +1\right)$. Thus, in this case we also have $\alpha>\frac{1}{p}$.

Suppose that $\nu \geq -\frac{1}{2}$. As it was mentioned above, the kernel $W_t^{\nu }(x,y)$, $x,y\in (0,1)$ and $t>0$, has a Gaussian upper bound. Then, according to \cite[Lemma 2.2]{DOS}
\begin{equation}\label{1.2}
\int_0^1\left|K_{F(\sqrt{L_\nu})}(x,y)\right|^2dy\leq \frac{C}{|I(x,t^{-1})|}\|F_t\|_\infty ^2,
\end{equation}
for every $t>0$ and every Borel function $F$ on $\mathbb{R}$ such that $\supp (F)\subset [0,t]$. Here, for each $x\in (0,1)$ and $t>0$, $I(x,t)=(0,1)\cap (x-t,x+t)$, and $|I(x,t)|$ is the length of the interval $I(x,t)$. Also, for every Borel function $F$ on $\mathbb{R}$, $K_{F(\sqrt{L_\nu })}$ represents the integral kernel of the operator $F(\sqrt{L_\nu})$ that is defined by functional calculus, and $F_t(z)=F(tz)$, for $z\in \mathbb{R}$ and $t\in (0,\infty )$.

In the above inequality \eqref{1.2}, $\|\cdot\|_\infty$ cannot be replaced by $\|\cdot\|_p$ for any $p<\infty$, since, for every $j\in \mathbb{N}$, $s_{\nu ,j}^2$ belongs to the pointwise spectrum of $L_\nu$ (i.e. the spectrum of $L_\nu$ is non-empty; see \cite[p. 450]{DOS}).

From \cite[Theorem 1.1]{CDY} we deduce that \eqref{1.1} holds when $1<p<\infty$ and $\alpha >2\left|\frac{1}{p}-\frac{1}{2}\right|+\frac{1}{2}$. Note that 
$$
2\left|\frac{1}{p}-\frac12 \right|+\frac12 =
\left\{\begin{array}{ll}
        \frac{2}{p}-\frac12 ,&{ p\in (1,2)},\\
        \\
        \frac{3}{2}-\frac{2}{p},&{ p\in [2,\infty)},
        \end{array}
    \right.
$$
and so, $2\left|\frac{1}{p}-\frac12 \right|+\frac12 \geq \max\left\{\frac{1}{p},1-\frac{1}{p}\right\}$. 

On the other hand, if $-1<\nu <-\frac{1}{2}$, then $W_t^{L_\nu }(x,y)$, $x,y\in (0,1)$ and $t>0$, has not got Gaussian upper bounds. Hence, Theorem \ref{Theorem1.1} generalizes \cite[Theorem 1.1]{CDY} for the operator $L_\nu$. The proof of Theorem \ref{Theorem1.1}  is different from \cite[Theorem 1.1]{CDY}, since we take advantage of the properties of several Bessel--type functions.

One of the steps in the proof of Theorem \ref{Theorem1.1} consists in proving $L^p$--boundedness properties for the Stein's square functions defined through the Hankel transform instead of using the Fourier transform.

Let again $\nu >-1$. The Hankel transform $h_\nu (f)$ of a complex measurable function $f$ on $(0, \infty )$ such that $x^{\min\{\nu +\frac{1}{2}, 0\}}f\in L^1(0,1)$ and $f\in L^1(1,\infty )$ is defined by
$$
h_\nu (f)(x)=\int_0^\infty \sqrt{xy}J_\nu (xy)f(y)dy,\quad x\in (0,\infty ).
$$

The transform $h_\nu$ can be extended from $C_c^\infty (0,\infty)$, the space of smooth functions with compact support in $(0,\infty )$, to $L^2(0,\infty )$ as an isometry on $L^2(0,\infty )$ (\cite[Lemma 2.7]{BS}). We define the Bessel operator $B_\nu$ as follows
$$
B_\nu f=h_\nu \left(x^2h_\nu (f)\right),\quad f\in D(B_\nu),
$$
where the domain $D(B_\nu)$ of $B_\nu$ is given by
$$
D(B_\nu)=\left\{f\in L^2(0,\infty ): x^2h_\nu (f)\in L^2(0,\infty )\right\}.
$$
Since $C_c^\infty (0,\infty )\subset D(B_\nu)$, $D(B_\nu)$ is a dense subspace of $L^2(0,\infty )$. If $f\in C_c^\infty (0,\infty )$, we can see that 
$$
B_\nu f=\mathbb{L}_\nu f=-x^{-\nu -\frac{1}{2}}\frac{d}{dx}x^{2\nu +1}\frac{d}{dx}x^{-\nu -\frac{1}{2}}f.
$$
We have that $B_\nu $ is a positive self-adjoint operator and its spectrum is $\sigma (B_\nu)=[0,\infty )$. Hence, $B_\nu$ is a sectorial operator and $-B_\nu$ generates a holomorphic semigroup of operators $\{W_t^{B_\nu }\}_{t>0}$. For every $t>0$ and $f\in L^2(0,\infty )$ we can write
$$
W_t^{B_\nu}(f)(x)=\int_0^\infty W_t^{B_\nu }(x,y)f(y)dy,\quad x\in (0,\infty ),
$$
where
$$
W_t^{B_\nu }(x,y)=\frac{\sqrt{xy}}{2t}I_\nu \left(\frac{xy}{2t}\right)e^{-\frac{x^2+y^2}{4t}},\quad x,y,t\in (0,\infty ).
$$
Here, $I_\nu$ denotes the modified Bessel function of the first kind and order $\nu$.

If $\nu \geq -\frac{1}{2}$, $W_t^{B_\nu }(x,y,t)$, for $x,y,t\in (0,\infty )$, has a Gaussian upper bound, but, if $-1<\nu <-\frac{1}{2}$, it lacks this property.

We define the Bochner--Riesz means associated with $B_\nu$ of order $\alpha >0$ by
$$
R_t^{\alpha ,B_\nu }f=h_\nu \left(\left(1-\frac{z^2}{t^2}\right)_+^\alpha h_\nu (f)\right),\quad f\in L^2(0,\infty ).
$$
The Stein's square function for $B_\nu$ of order $\alpha$ is defined by
$$
G_\alpha ^{B_\nu}(f)(x)=\left(\int_0^\infty \left|\frac{\partial }{\partial _t}R_t^{\alpha ,B_\nu }f(x)\right|^2tdt\right)^{\frac{1}{2}}.
$$
If $\nu>-\frac{1}{2}$ we can prove that
$$
\int_0^\infty |K_{F(\sqrt{B_\nu})}(x,y)|^2dy\leq \frac{C}{|J(x,t^{-1})|}\|F_t\|_2^2,
$$
for every $x,t>0$ and every Borel function $F$ on $\mathbb{R}$ such that $\supp (F)\subset [0,t]$, where $K_{F(\sqrt{B_\nu})}$ represents the integral kernel of the operator $F(\sqrt{B_\nu})$ that is defined by functional calculus and, as before, $F_t(z)=F(tz)$, for $z\in \mathbb{R}$ and $t>0$. Here, for each $x,t\in (0,\infty )$, $|J(x,t)|$ is the length of the interval $J(x,t)=(x-t,x+t)\cap (0,\infty )$.

According to \cite[Theorem 1.1]{CDY}, $G_\alpha ^{B\nu}$ is bounded on $L^p(0,\infty )$ when $\nu >-\frac{1}{2}$ and 
$$
\alpha >\left|\frac{1}{p}-\frac12 \right|+\frac12 =
\left\{\begin{array}{ll}
        \displaystyle \frac{1}{p},&p\in (1,2),\\
        \displaystyle 1-\frac{1}{p},&p\in [2,\infty ).
        \end{array}
\right.
$$
We improve this result as follows.

\begin{thm}\label{Theorem1.2}
$G_\alpha ^{B_\nu}$ is bounded on $L^p(0,\infty )$, when $\alpha >\max\left\{\frac{1}{p},\frac12 \right\}$ and
\begin{itemize}
    \item[(i)] $\nu >-\frac12 $ and $1<p<\infty$ or
    \item[(ii)]$-1<\nu \leq -\frac12 $ and $\frac{2}{2\nu +3}<p<-\frac{2}{2\nu +1}$.
\end{itemize}
Moreover, if $G_\alpha ^{B_\nu}$ is bounded on $L^p(0,\infty )$ with $\nu >-1$, $\frac{2}{p}>-(2\nu +1)$ and $1<p<\infty$, then $\alpha >\frac{1}{p}$.
\end{thm}

In order to prove this result we use a vector--valued version of a transplantation result for Hankel transforms due to Nowak and Stempak (\cite[Theorem 2.1]{NS}).

In \cite{BS} some transference results of $L^p$--boundedness properties from Fourier--Bessel to Hankel multipliers were established. We now show a transference of $L^p$--boundedness from $G_\alpha ^{L_\nu}$ to $G_\alpha ^{B_\nu}$. In order to do so, we consider square functions as vector--valued multipliers. The procedure developed in the proof of our transference property between square functions can also be used to see, in a simpler way, the transference results for Littlewood-Paley squares functions stated in \cite{CFS} and \cite{NU}.
\begin{thm}\label{Theorem1.3}
Let $\nu >-\frac{1}{2}$, $1<p<\infty$ and $\alpha >\frac{1}{2}$. If $G_\alpha ^{L_\nu}$ is bounded on $L^p(0,1)$, then $G_\alpha ^{B_\nu}$ is bounded on $L^p(0,\infty )$. In particular, if $1<p\leq 2$ and $G_\alpha ^{L_\nu}$ is bounded on $L^p(0,1)$, then $\alpha >\frac{1}{p}$.
\end{thm}

By combining the above theorems it follows that when $\nu>-\frac{1}{2}$ and $1<p\leq 2$, $G_\alpha ^{L_\nu}$ (respectively, $G_\alpha ^{B_\nu}$) is bounded on $L^p(0,1)$  (respectively, $L^p(0,\infty )$) if and only if $\alpha >\frac{1}{p}$. We recall that the classical Stein's square function $G_\alpha$ is bounded on $L^p(\mathbb{R}^n)$, with $1\leq p\leq 2$, if and only if $\alpha >n\left(\frac{1}{p}-\frac{1}{2}\right)+\frac{1}{2}$ (see \cite{LRS1} and \cite{Su1}).

Let $\nu >-1$. If $m=\{m_j\}_{j=1}^\infty$ is a complex sequence, the spectral multiplier associated with $L_\nu$ defined by $m$ is given by
$$
m(L_\nu)f=\sum_{j=1}^\infty m_jc_{\nu ,j}(f)\phi _{\nu ,j},\quad f\in L^2(0,1).
$$
In the case where $m$ is a bounded sequence, $m(L_\nu)$ is a bounded operator on $L^2(0,1)$. In the following theorem we give a smoothness condition on $m$ that allows us to extend $m(L_\nu)$ from $L^2(0,1)\cap L^p(0,1)$ to $L^p(0,1)$ as a bounded operator on $L^p(0,1)$, with $1<p<\infty$. The Stein's square functions $G_\alpha ^{L_\nu}$ plays a fundamental role in its proof.

\begin{thm}\label{Theorem1.4}
Suppose that $m=\{m_j\}_{j=1}^\infty$ is a complex sequence such that
$$
|||m|||=\sup_{j\in \mathbb{N}}|m_j|+\left(\sum_{\ell =1}^\infty \frac{1}{\ell ^s}\sum_{k=1}^\ell \left|m(s_{\nu ,k+1})-m(s_{\nu ,k})\right|^2\right)^{\frac{1}{2}}<\infty.
$$
If $1<p<\infty$ and some of the following conditions hold
\begin{itemize}
\item[(i)]$\nu >-\frac12 $;
\item[(ii)]$-1<\nu \leq -\frac12 $ and $\frac{2}{2\nu +3}<p\leq 2$;
\item[(iii)] $-1<\nu \leq -\frac12 $ and $2\leq p<-\frac{1}{\nu +\frac{1}{2}}$;
\end{itemize}
the operator $m(L_\nu)$ can be extended from $L^2(0,1)\cap L^p(0,1)$ to $L^p(0,1)$ as a bounded operator on $L^p(0,1)$.
\end{thm}

In the following sections we present proofs of our results. Throughout this paper $C$ denotes a positive constant that can change from one line to another.

\section{Proof of Theorem \ref{Theorem1.2}}
In order to prove Theorem \ref{Theorem1.2}, let $\nu >-1$ and recall that, for every $t>0$,
$$
R_t^{\alpha ,B_\nu}f=h_\nu \left(\left(1-\frac{z^2}{t^2}\right)_+^\alpha h_\nu (f)\right),\quad f\in L^2(0,\infty ).
$$
Since $\sqrt{z}J_{-\frac{1}{2}}(z)=\sqrt{2/\pi}\cos z$, $z\in (0,\infty )$, we have that $h_{-\frac{1}{2}}=\mathcal{F}_{\rm c}$, the Fourier-cosine transformation, given by
$$
\mathcal{F}_{\rm c}(f)(x)=\sqrt{\frac{2}{\pi }}\int_0^\infty \cos (xy)f(y)dy,\quad x\in (0,\infty ),
$$
for every $f\in L^1(0,\infty )$. It is well-known that this transform can be extended from $L^1(0,\infty)\cap L^2(0,\infty)$ as a bounded operator on $L^2(0,\infty)$. 
We get, for each $t>0$, that
$$
R_t^{\alpha ,B_{-\frac{1}{2}}}f=\mathcal{F}_{\rm c}\left(\left(1-\frac{z^2}{t^2}\right)_+^\alpha \mathcal{F}_{\rm c}(f)\right),\quad f\in L^2(0,\infty ).
$$

We consider the Fourier transformation defined, as usual, by
$$
\widehat{f}(x)=\frac{1}{\sqrt{2\pi }}\int_{-\infty }^{+\infty }e^{-ixy}f(y)dy,\quad x\in \mathbb{R},
$$
for every $f\in L^1(\mathbb{R})$. It is also widely known that this transformation can be extended from $L^1(\mathbb{R})\cap L^2(\mathbb{R})$ to $L^2(\mathbb{R})$ as an isometry on $L^2(\mathbb{R})$.

If $g$ is a measurable function defined on $(0,\infty )$ we consider $g_{\rm e}$ as the even extension of $g$ to $\mathbb{R}$. For every $f\in L^2(0,\infty )$, we have that,
$$
(\mathcal{F}_{\rm c}(f))_{\rm e}=\widehat{f_{\rm e}},
$$
and it follows that
$$
R_t^{\alpha ,B_{-\frac{1}{2}}}f=\left(\left(1-\frac{z^2}{t^2}\right)_+^\alpha\widehat{f_{\rm e}}\right)^{\widecheck{\;}}=R_t^\alpha f_{\rm e}\bigg|_{(0,\infty )}, \quad t>0.
$$
Note that $\widehat{f_{\rm e}}$ is an even function and $\widecheck{g}=(\mathcal{F}_{\rm c}(g))_{\rm e}$, provided that $g\in L^2(\mathbb{R})$ is even.

We can write 
$$
R_t^{\alpha ,B_\nu}f=T_{\nu ,-\frac{1}{2}}R_t^{\alpha ,B_{-\frac{1}{2}}}T_{-\frac{1}{2},\nu }f,\quad f\in L^2(0,\infty ),
$$
where $T_{\nu ,-\frac{1}{2}}$ and $T_{-\frac{1}{2},\nu}$ represent the following transplantation operators
$$
T_{\nu ,-\frac{1}{2}}=h_\nu h_{-\frac{1}{2}}\quad \mbox{and} \quad T_{-\frac{1}{2},\nu }=h_{-\frac{1}{2}}h_\nu .
$$
To simplify, we put $\mathcal{H}=L^2((0,\infty ),\frac{dt}{t})$. We consider
$\widetilde{G}_\alpha^{B_{-\frac{1}{2}}}: L^2(0,\infty )\longrightarrow L^2((0,\infty ),\mathcal{H})$ defined by
$$
\widetilde{G}_\alpha ^{B_{-\frac{1}{2}}}(f)(x,t)=t\frac{\partial }{\partial t}R_t^{\alpha ,B_{-\frac{1}{2}}}f(x),\quad f\in L^2(0,\infty )\mbox{ and }x,t\in (0,\infty ).
$$
According to \cite{Su2, Su1, Su3}, $\widetilde{G}_\alpha ^{B_{-\frac{1}{2}}}$ is bounded from $L^p(0,\infty )$ into $L^p((0,\infty ),\mathcal{H})$ whenever $p\in (1,\infty)$ and $\alpha >\max\left\{\frac{1}{2},\frac{1}{p}\right\}$.

We also consider, for every $\mu >-1$, the operator $\widetilde{G}_\alpha^{B_\mu}: L^2(0,\infty )\longrightarrow L^2((0,\infty ),\mathcal{H})$, defined by
$$
\widetilde{G}_\alpha ^{B_\mu}(f)(x,t)=t\frac{\partial }{\partial t}R_t^{\alpha ,B_\mu }f(x),\quad f\in L^2(0,\infty )\mbox{ and }x,t\in (0,\infty ).
$$
If $\alpha >0$ and $1<p<\infty$, $G_\alpha^{B_\nu}$ is bounded on $L^p(0,\infty)$ if and only if $\widetilde{G}_\alpha ^{B_\nu}$ is bounded from $L^p(0,\infty )$ into $L^p((0,\infty), \mathcal{H})$.

According to \cite[Theorem 2.1]{NS} the operator $T_{-\frac{1}{2},\nu}$ is bounded on $L^p (0,\infty)$ provided that some of the following conditions hold:

\begin{enumerate}[label=\emph{(\roman*)}]
\item $\nu +\frac12 =2k$, for some $k\in \mathbb{N} \cup \{0\}$, and $1<p<\infty$;
\item $\nu +\frac12 \not=2k$, for every $k\in \mathbb{N} \cup \{0\}$, and  $\frac{2}{2\nu +3}<p<\infty$.
\end{enumerate}
\noindent Note that, whatever the case, if $\nu \geq -\frac12 $, then $\frac{2}{2\nu +3}<p<\infty$.

On the other hand, the operator $T_{\nu ,-\frac{1}{2}}$ is bounded on $L^2((0,\infty ),\mathcal{H})$. Indeed, for every $\mu >-1$, the Hankel transformation $h_\mu$ is an isometry in $L^2(0,\infty)$ (\cite[Lemma 2.7]{BS}). Then, $h_\mu$ can be extended from $L^2(0,\infty)\otimes \mathcal H$ to $L^2((0,\infty ),\mathcal{H})$ as a bounded operator on $L^2((0,\infty ),\mathcal{H})$ (see, for instance, \cite[Theorem 1.8.2]{ABHN}).

By proceeding as in the proof of \cite[Theorem 2.1]{NS} but in the $\mathcal{H}$--valued setting we conclude that $T_{\nu ,-\frac{1}{2}}$ is bounded on $L^p((0,\infty),\mathcal{H})$, $1<p<\infty$, when some of the following properties are satisfied:
\begin{enumerate}[label=\emph{(\roman*)}]
\item $\nu +\frac12 =2k$, for some $k\in \mathbb{N} \cup \{0\}$;
\item $\nu +\frac12 \not=2k$, for every $k\in \mathbb{N} \cup \{0\}$ and  $\frac{2}{p}>-(2\nu +1)$.
\end{enumerate}
\noindent Note that if $\nu >-\frac12 $, then $\frac{2}{p}>-(2\nu +1)$, for every $1<p<\infty$.

By combining the above results we conclude that $G_\alpha ^{B_\nu}$ is bounded on $L^p(0,\infty )$ when
\begin{enumerate}[label=\emph{(\roman*)}]
\item $\nu >-\frac12 $, provided that $1<p<\infty$ and $\alpha >\max\left\{\frac12 ,\frac{1}{p}\right\}$; or when
\item $-1<\nu \leq -\frac12 $, $\frac{2}{2\nu +3}<p<-\frac{2}{2\nu +1}$ and $\alpha >\max\left\{\frac12 , \frac{1}{p}\right\}$.
\end{enumerate}
Thus, the proof of the first part of Theorem \ref{Theorem1.2} is finished.

In order to prove the necessity of the conditions, we consider a function $\phi\in C^\infty (0,\infty )$ such that
$\phi (x)=1$ for $x\in (0,1)$, and $\phi (x)=0$ for $x\in (2,\infty )$. Let $\nu >-1$. It is clear that $x^{\nu+\frac{1}{2}}\phi \in L^1(0,\infty )$ and $\phi \in L^1(0,\infty)\cap L^2(0,\infty)$. We define $\psi=h_\nu (\phi )$. We claim that
\begin{equation}\label{psi}
|\psi (x)|\leq C\left\{\begin{array}{ll}
                x^{\nu +\frac{1}{2}},&0<x\leq 2,\\
                &\\
                \displaystyle \frac{1}{x},&x>2.
                \end{array}
\right.
\end{equation}
Indeed, since $z^{-\nu }J_\nu (z)$ and $\sqrt{z}J_\nu (z)$ are bounded functions on $(0,1)$ and on $(1,\infty)$, respectively (\cite[pp. 104 and 122, respectively]{Le}), we get
$$
|\psi (x)|\leq C\int_0^2(xy)^{\nu +\frac{1}{2}}dy\leq Cx^{\nu +\frac{1}{2}},\quad 0<x\leq 2,
$$
and
$$
\left|\int_0^{\frac{1}{x}}\sqrt{xy}J_\nu (xy)\phi (y)dy\right|\leq \frac{C}{x},\quad x>2.
$$
Also, by taking into account that $\frac{d}{dx}[x^{-\nu }J_\nu (x)]=-x^{-\nu }J_{\nu +1}(x)$ (\cite[(5.3.5)]{Le}) and $\mathbb{L}_\nu (\sqrt{x}J_\nu (x))=\sqrt{x}J_\nu (x)$, $x\in (0,\infty )$, (\cite[(5.3.7)]{Le}), partial integration leads to
\begin{align*}
\int_a^b\sqrt{xy}J_\nu (xy)\phi (y)dy&=-\frac{1}{x^2}\int_a^by^{-\nu -\frac{1}{2}}\frac{d}{dy}\left[y^{2\nu +1}\frac{d}{dy}\left(y^{-\nu -\frac{1}{2}}\sqrt{xy}J_\nu (xy)\right)\right]\phi (y)dy\\
&=-\frac{1}{x^2}\left(-x\sqrt{xy}J_{\nu +1}(xy)\phi (y)\bigg\rvert_a^b-y^{\nu +\frac{1}{2}}\sqrt{xy}J_\nu (xy)\frac{d}{dy}\left(y^{-\nu -\frac{1}{2}}\phi (y)\right)\bigg\rvert_a^b\right.\\
&\quad + \left.\int_a^b\sqrt{xy}J_\nu (xy)y^{-\nu -\frac{1}{2}}\frac{d}{dy}\left(y^{2\nu +1}\frac{d}{dy}\left(y^{-\nu -\frac{1}{2}}\phi (y)\right)\right)dy\right), 
\end{align*}
for any $0<a\leq b<\infty$ and $x\in (0,\infty )$. It yields 
\[
\left|\int_{\frac{1}{x}}^2\sqrt{xy}J_\nu (xy)\phi (y)dy\right|\leq C\frac{x+1}{x^2}\leq \frac{C}{x},\quad x>2.
\]
We conclude that $|\psi (x)|\leq \frac{C}{x}$, $x>2$, so \eqref{psi} holds. 
This inequality let us deduce that $\psi \in L^p(0,\infty )$, provided that $1<p<\infty $ and $\nu +\frac12  >-\frac{1}{p}$.

We can write 
\begin{align*}
\widetilde{G}_\alpha ^{B_\nu }(\psi )(x,t)&=2\alpha h_\nu \left(\left(1-\frac{z^2}{t^2}\right)_+^{\alpha -1}\frac{z^2}{t^2}h_\nu (\psi )(z)\right)(x)=2\alpha h_\nu \left(\left(1-\frac{z^2}{t^2}\right)_+^{\alpha -1}\frac{z^2}{t^2}\phi (z)\right)(x)\\
&=2 h_\nu \left(\left(1-\frac{z^2}{t^2}\right)_+^{\alpha -1}\frac{z^2}{t^2}\right)(x),\quad x\in (0,\infty),\;t\in (0,1).
\end{align*}
According to \cite[(34) p. 26]{EMOT}, for $t\in (0,1)$ and $x\in(0,\infty )$, we get 
\begin{align*}
\widetilde{G}_\alpha ^{B_\nu }(\psi )(x,t)&=\frac{2\alpha }{t^{2\alpha }}\int_0^t\left(t^2-z^2\right)^{\alpha -1}z^2\sqrt{xz}J_\nu (xz)dz\\
&=\frac{2\alpha B(\alpha , \frac{7}{4}+\frac{\nu}{2})}{2^{\nu +1}\Gamma (\nu +1)}t^{\nu +\frac{3}{2}}x^{\nu +\frac12 }\;_1F_2\left(\frac{7}{4}+\frac{\nu}{2}; \nu +1,\alpha +\frac{7}{4}+\frac{\nu}{2}; -\frac{x^2t^2}{4}\right).
\end{align*}
Here $B$ represents the Beta Eulerian function and $\!\!\;_1F_2$ denotes the hypergeometric function.

By \cite{WR}, for certain $a,b\in \mathbb{R}$ we have that
\begin{align*}
\widetilde{G}_\alpha ^{B_\nu }(\psi )(x,t)&=at^{-2}x^{-3}\left(1+O\left(\frac{1}{x^2t^2}\right)\right)\\
&\quad +bt^{-\alpha +1}x^{-\alpha }\left(\cos \left(xt-\frac{\pi}{2}\left(\nu+\frac12 +\alpha \right)\right)+O\left(\frac{1}{xt}\right)\right),\quad t\in \left(\tfrac{1}{2},1\right), \;x\in [1,\infty ).
\end{align*}
Thus, it follows that 
$$
\int_1^\infty \left(\int_{\frac{1}{2}}^1 \left|\widetilde{G}_\alpha ^{B_\nu }(\psi )(x)\right|^2\frac{dt}{t}\right)^{p/2}dx<\infty,
$$
if and only if
$$
\int_1^\infty \left(\int_{\frac{1}{2}}^1\cos ^2\left(xt-\frac{\pi}{2}(\nu+\frac12 +\alpha )\right)dt\right)^{p/2}x^{-\alpha p}dx<\infty.
$$
We can write
\begin{align*}
\int_{\frac{1}{2}}^1\cos ^2\left(xt-\frac{\pi}{2}\left(\nu+\frac12 +\alpha\right)\right)dt\\
&\hspace{-3cm}=\frac{1}{4}\left[1+\frac{1}{x}\left(\sin \left(2x-\pi \left(\nu +\frac {1}{2}+\alpha \right)\right)-\sin \left(x-\pi \left(\nu +\frac12 +\alpha\right)\right)\right)\right],\quad x\in [1,\infty ).
\end{align*}
Then, there exists $x_0>1$ such that
$$
\int_{\frac{1}{2}}^1\cos ^2\left(xt-\frac{\pi}{2}\left(\nu+\frac12 +\alpha \right)\right)dt\geq \frac{1}{8},\quad x>x_0.
$$
We conclude that
$$
\int_1^\infty \left(\int_{\frac{1}{2}}^1\cos ^2\left(xt-\frac{\pi}{2}\left(\nu+\frac12 +\alpha \right)\right)dt\right)^{p/2}x^{-\alpha p}dx<\infty
$$
only if $\alpha >\frac{1}{p}$.

Hence, if $1<p<\infty$, $\nu +\frac{1}{2}>-\frac{1}{p}$ and $\widetilde{G}_\alpha ^{B_\nu }$ is bounded from $L^p(0,\infty )$ into itself, necessarily $\alpha >\frac{1}{p}$. Thus, the proof is complete.

\section{Proof of Theorem \ref{Theorem1.1}}

Our first objective is to prove that, if any of the conditions $(i)$, $(ii)$ or $(iii)$ is satisfied, then
\begin{equation}\label{Objective1}
\|G_\alpha ^{L_\nu }(f)\|_p\leq C\|f\|_p,\quad f\in L^2(0,1)\cap L^p(0,1).
\end{equation}
Assume that $\nu >-1$ and $f\in L^2(0,1)$. We can write
$$
R_t^{\alpha ,L_\nu}f(x)=\int_0^1R_t^{\alpha,L_\nu}(x,y)f(y)dy,\quad x\in (0,1),\; t>0,
$$
where
$$
R_t^{\alpha , L_\nu}(x,y)=\sum_{j=1}^\infty \left(1-\frac{s_{\nu ,j}^2}{t^2}\right)_+^{\alpha}\phi _{\nu ,j}(x)\phi _{\nu ,j}(y),\quad x,y\in (0,1), \;t>0.
$$

It is clear that
\begin{align*}
t\frac{\partial}{\partial t}R_t^{\alpha ,L_\nu }(x,y)&=2\alpha \sum_{j=1}^\infty \left(1-\frac{s_{\nu ,j}^2}{t^2}\right)_+^{\alpha -1}\frac{s_{\nu ,j}^2}{t^2}\phi _{\nu ,j}(x)\phi _{\nu ,j}(y)\\
&=2\alpha (R_t^{\alpha -1,L_\nu}(x,y)-R_t^{\alpha ,L_\nu }(x,y)),\quad x,y\in (0,1), \;t\in (0,\infty)\setminus\{s_{\nu ,j}\}_{j=1}^\infty.
\end{align*}

By proceeding as in the proof of \cite[Lemma 1]{CR1} we can see that, for every $x,y\in (0,1)$ and $t>0$,
$$
t\frac{\partial}{\partial t}R_t^{\alpha ,L_\nu }(x,y)=H_{\nu ,t;1}^\alpha (x,y)+H_{\nu ,t;2}^\alpha (x,y),
$$
where
$$
H_{\nu ,t;1}^\alpha (x,y)=2\alpha \sqrt{xy} \int_0^tz \left(1-\frac{z^2}{t^2}\right)^{\alpha -1}\left(\frac{z}{t}\right)^2J_\nu (xz)J_\nu (yz)dz,
$$
and
\begin{equation}\label{Ht2}
H_{\nu ,t;2}^\alpha (x,y)=2\alpha \lim_{\varepsilon \rightarrow 0^+}\frac{\sqrt{xy}}{2}\int_{S_{\varepsilon,t}} \left(1-\frac{z^2}{t^2}\right)^{\alpha -1}\left(\frac{z}{t}\right)^2\frac{zH_\nu ^{(1)}(z)J_\nu (xz)J_\nu (yz)}{J_\nu (z)}dz.
\end{equation}
Here, $S_{\varepsilon, t}$, $\varepsilon , t>0$, denotes the path defined as follows 
$$
S_{\varepsilon, t}(u)=\left\{\begin{array}{lll}
t+iu, & & u\in(\varepsilon,+\infty), \\
 & & \\
-t-iu, & & u\in(-\infty,-\varepsilon), 
\end{array}\right.
$$
and $H_\nu ^{(1)}$ is the Hankel function of the first kind and order $\nu$.

A few words about the definitions of the functions in the complex plane. We need to precise in every case the branch of the logarithm that is considered. The Bessel function $J_\nu$ is defined, when $\nu\in \mathbb{R}$, by
$$
J_\nu (z)=\sigma_\nu (z)\sum_{k=0}^\infty \frac{(-1)^k2^{-\nu -2k}z^{2k}}{\Gamma (k+1)\Gamma (k+\nu +1)},\quad z\in \mathbb{C}\setminus i(-\infty,0],
$$
where $\sigma_\nu(z)=e^{\nu (\ln|z|+i\theta (z))}$, $z\in \mathbb{C}\setminus i(-\infty,0]$, and
$$
\theta (z)=\left\{
\begin{array}{ll}
{\rm Arg}\;z,&-\frac{\pi}{2}<{\rm Arg}\;z\leq \pi,\\
2\pi +\rm{Arg}\;z,&-\pi <{\rm Arg}\;z<-\frac{\pi}{2}.
\end{array}
\right.
$$
The function $H_\nu ^{(1)}$ is defined, for every $\nu >-\frac{1}{2}$, $\nu \not \in\mathbb{N}$, by
$$
H_\nu ^{(1)}(z)=\frac{J_{-\nu }(z)-e^{-\nu \pi i}J_\nu (z)}{i\sin (\nu \pi)},\quad z\in \mathbb{C}\setminus i(-\infty ,0].
$$
When $\nu =n\in \mathbb{N}$, $H_\nu ^{(1)}(z)=\lim_{\nu \rightarrow n}H_\nu ^{(1)}(z)$, $z\in \mathbb{C}\setminus i(-\infty ,0]$.

Finally, for every $\beta \in \mathbb{R}\setminus\mathbb{Z}$, the function $w_\beta(z)=(1-z^2)^\beta$ is defined in a holomorphic way in $\mathbb{C}\setminus ((-\infty,-1]\cup [1,+\infty))$.

Thus, the function
$$
\xi (z)=\left(1-\frac{z^2}{t^2}\right)^{\alpha -1}\frac{z^2}{t^2}\frac{zH_\nu ^{(1)}(z)J_\nu (xz)J_\nu (yz)}{J_\nu (z)},
$$
is holomorphic in $\mathbb{C}\setminus ((-\infty,-t]\cup [t,+\infty)\cup (i(-\infty,0])\cup \{s_{\nu ,j}\}_{j=1}^\infty )$. $\xi$ has branch points in $=-t$, $z=t$, and $z=0$ and $\xi$ has poles when $z=\pm s_{\nu ,j}$, $j=1,...,j(t)$, where $j(t)\in \mathbb{N}$ and $s_{\nu ,j(t)+1}\geq t$.

We choose an even and smooth function $\phi$ on $\mathbb{R}$ such that, $0\leq \phi \leq 1$, $\phi (x)=1$, for $x\in [-1,1]$, and $\phi (x)=0$ for $x\not\in (-2,2)$, and write
\begin{align}\label{decomposition}
t\frac{\partial}{\partial t}R_t^{\alpha ,L_\nu }(x,y)&=H_{\nu,t;1}^\alpha (x,y)+\left(1-\phi (t|x-y|)\right)H_{\nu,t;2}^\alpha (x,y)+\phi (t|x-y|)H_{\nu,t;2}^\alpha (x,y)\nonumber\\
&=:\sum_{j=1}^3A_{\nu ,t;j}^\alpha (x,y), \quad x,y\in (0,1),\;t>0.
\end{align}

For $j=1,2,3$, we consider the operators
$$
\mathbb{A}_{\nu ,t;j}^\alpha (f)(x)=\int_0^1A_{\nu ,t;j}^\alpha (x,y)f(y)dy,\quad x\in (0,1),\;t>0,
$$
and the square functions 
$$
G\mathbb{A}_{\nu ;j}^\alpha (f)(x)=\left(\int_0^\infty |\mathbb{A}_{\nu ,t;j}^\alpha (f)(x)|^2\frac{dt}{t}\right)^{\frac{1}{2}},\quad x\in (0,1).
$$

In order to show \eqref{Objective1} we will analyze first the $L^p$--boundedness properties for $G\mathbb{A}_{\nu ;j}^\alpha $, $j=1,2,3$.

\noindent $\bullet$ {\bf $L^p$--boundedness properties for $G\mathbb{A}_{\nu ;1}^\alpha$}. 
By defining $\widetilde{f}$ as the extension of $f$ to $(0,\infty)$ such that $\widetilde{f}(x)=0$ for $x\in (1,\infty )$, we get
$$
G\mathbb{A}_{\nu ;1}^\alpha (f)(x)\leq CG_\alpha ^{B_\nu}(\widetilde{f})(x),\quad x\in (0,1).
$$
Then, according to Theorem \ref{Theorem1.2}, $G\mathbb{A}_{\nu;1}^\alpha$ is bounded on $L^p(0,1)$ when $\alpha >\max\left\{\frac{1}{p},\frac12 \right\}$ and one of the following conditions is satisfied:
\begin{itemize}
\item[($i'$)] $\nu \geq -\frac12 $ and $1<p<\infty$;

\item[($ii'$)] $-1<\nu < -\frac12 $ and $\frac{2}{2\nu +3}<p<-\frac{2}{2\nu +1}$.
\end{itemize}
\vspace*{0.3cm}

\noindent $\bullet$ {\bf $L^p$--boundedness properties for $G\mathbb{A}_{\nu;2}^\alpha$}. We are going to prove that $G\mathbb{A}_{\nu;2}^\alpha$ is a bounded operator on $L^p(0,1)$, when $\alpha >1$ and one of the above conditions $(i')$ or $(ii')$ hold.

For every $\varepsilon, t>0$, let us represent by $S_{\varepsilon,t}^+$ and $S_{\varepsilon ,t}^-$ the paths given by
$$
S_{\varepsilon,t}^+(u)=t+iu, \quad u\in (\varepsilon ,+\infty)\quad \mbox{ and} \quad S_{\varepsilon, t}^-(u)=-t-iu, \quad u\in (-\infty,-\varepsilon).
$$
It is clear that $z\in S_{\varepsilon,t}^-$ if and only if $-\overline{z}\in S_{\varepsilon, t}^+$. Moreover, by using that $J_\nu (-z)=e^{\nu \pi i}J_\nu (z)$, $H_\nu ^{(1)}(-z)=-e^{-\nu \pi i}H_\nu ^{(2)}(z)$, $J_\nu (\overline{z})=\overline{J_\nu (z)}$ and $H_\nu ^{(2)}(\overline{z})=\overline{H_\nu ^{(1)}(z)}$ for $z\in \mathbb{C}\setminus (-\infty ,0]$, where $H_\nu ^{(2)}$ denotes the Hankel function of the second kind and order $\nu$, we have that, if $z\in S_{\varepsilon, t}^-$, then
$$
\frac{H_\nu ^{(1)}(z)}{J_\nu (z)}\sqrt{zx}J_\nu (zx)\sqrt{zy}J_\nu (zy)
=\overline{\frac{H_\nu ^{(1)}(-\overline{z})}{J_\nu (-\overline{z})}\sqrt{-\overline{z}x}J_\nu (-\overline{z}x)\sqrt{-\overline{z}y}J_\nu (-\overline{z}y)},\quad x,y \in (0,1).
$$

Thus, we can write
\begin{align*}
\int_{S_{\varepsilon,t}^-} \left(1-\frac{z^2}{t^2}\right)^{\alpha-1}\frac{z^2}{t^2}\frac{H_\nu^{(1)}(z)}{J_\nu(z)}\sqrt{zx}J_\nu(zx)\sqrt{zy}J_\nu(zy) dz&\\
&\hspace{-4cm}=\int_{S_{\varepsilon,t}^-} \overline{\left(1-\frac{(-\overline{z})^2}{t^2}\right)^{\alpha-1}\frac{(-\overline{z})^2}{t^2}\frac{H_\nu^{(1)}(-\overline{z})}{J_\nu(-\overline{z})}\sqrt{-\overline{z}x}J_\nu(-\overline{z}x)\sqrt{-\overline{z}y}J_\nu(-\overline{z}y)} dz\\
&\hspace{-4cm}=\int_{S_{\varepsilon,t}^+} \overline{\left(1-\frac{z^2}{t^2}\right)^{\alpha-1}\frac{z^2}{t^2}\frac{H_\nu^{(1)}(z)}{J_\nu(z)}\sqrt{zx}J_\nu(zx)\sqrt{zy}J_\nu(zy)}dz.
\end{align*}

By taking into account this equality, it is sufficient to analyze the boundedness on $L^p(0,1)$ for the operator $G\mathbb{A}_{\nu;2}^{\alpha, +}$ defined like the operator $G\mathbb{A}_{\nu;2}^\alpha$ but considering the kernel $H_{\nu,t;2}^{\nu ,+}$ as in \eqref{Ht2}, by replacing the path $S_{\varepsilon, t}$ by $S_{\varepsilon, t}^+$. As we will see, we need to use this trick because asymptotic expansions for Bessel type functions do not work uniformly on $\mathbb{C}\setminus(-\infty,0]$ (see \cite[\S 5.11]{Le} and \cite[p. 199]{Wat}).

As before, let us consider $\mathcal{H}=L^2\left((0,\infty),\frac{dt}{t}\right)$. Minkowski's integral inequality leads to
$$
G\mathbb{A}_{\nu;2}^{\alpha,+}(f)(x)\leq \int_0^1\Big\|(1-\phi ((\cdot)|x-y|))H_{\nu, \, \cdot\, ;2}^{\alpha ,+}(x,y)\Big\|_\mathcal{H}|f(y)|dy,\quad x\in (0,1).
$$

We claim that, for every $x,y\in (0,1)$, 
\begin{equation}\label{claim}
\Big\|(1-\phi ((\cdot)|x-y|))H_{\nu, \, \cdot\, ;2}^{\alpha ,+}(x,y)\Big\|_\mathcal{H}
\leq C\left\{\begin{array}{ll}
            \frac{1}{x},&\frac{x}{2}\leq y\leq2x,\\[0.3cm]
            \psi_\nu (x,y),&y<\frac{x}{2}\mbox{ or }y>2x,
            \end{array}
\right.
\end{equation}
where
$$
\psi_\nu (x,y)=
\left\{\begin{array}{ll}
        \frac{\min\{x,y\}^{\nu +\frac{1}{2}}}{\max\{x,y\}^{\nu +3/2}},& -1<\nu <-\frac12 ,\\[0.3cm]
        \frac{1}{\max\{x,y\}}, &\nu \geq -\frac12 .
\end{array}
\right.
$$
To establish these estimations we are going to see that there exists $C>0$ such that, for every $x,y\in (0,1)$ and $0<\varepsilon <1$,
$$
\int_{\frac{1}{|x-y|}}^\infty\left|\int_{S_{\varepsilon,t}^+}\left(1-\frac{z^2}{t^2}\right)^{\alpha -1}\frac{z^2}{t^2}\frac{H_\nu ^{(1)}(z)\sqrt{xz}J_\nu (xz)\sqrt{yz}J_\nu (yz)}{J_\nu (z)}dz\right|^2\frac{dt}{t}
\leq C\left\{\begin{array}{ll}
            \frac{1}{x^2},&\frac{x}{2}\leq y\leq 2x,\\[0.3cm]
            \psi_\nu (x,y)^2,&y<\frac{x}{2}\mbox{ or }y>2x.
            \end{array}
\right.
$$

Let $0<\varepsilon <1$ and consider $x,y\in (0,1)$. According to \cite[p. 199]{Wat} we have that
$$
\sqrt{w}J_\nu (w)=a_0\cos (w-\theta_\nu)+b_0\frac{\sin (w-\theta_\nu)}{w}+K(w),\quad w\in \mathbb{C}\setminus(-\infty ,0],
$$
where $\theta_\nu=\frac{\nu \pi}{2}+\frac{\pi}{4}$ and, for every $0<r <\pi$, there exists $C>0$ such that $|K(w)|\leq Ce^{|{\rm Im }\;w|}/|w|^2$, $|{\rm arg}\;w|<\pi -r$. Then, for each $z\in S_{\varepsilon,t}^+$, $t>0$, we can write 
\begin{align*}
\sqrt{xz}J_\nu (xz)\sqrt{yz}J_\nu (yz)&=c_1\cos((x+y)z-2\theta_\nu)+c_2\cos((x-y)z)\\
&\quad +c_3\frac{x+y}{xyz}\sin((x+y)z-2\theta_\nu)+c_4\frac{x-y}{xyz}\sin((x-y)z)+K(x,y,z),
\end{align*}
for certain $c_j\in \mathbb{R}$, $j=1,...,4$, and where
\begin{align*}
|K(x,y,z)|&=e^{(x+y){\rm |Im }\;z|}\left(O\left(\frac{1}{xy|z|^2}\right)+O\left(\frac{1}{x^2|z|^2}\right)+O\left(\frac{1}{y^2|z|^2}\right)\right.\\
&\quad \left.\quad +O\left(\frac{1}{xy^2|z|^3}\right)+
O\left(\frac{1}{x^2y|z|^3}\right)+O\left(\frac{1}{x^2y^2|z|^4}\right)\right).
\end{align*}
Thus, for every $t>0$, 
\begin{align*}
    \int_{S_{\varepsilon,t}^+} \left(1-\frac{z^2}{t^2}\right)^{\alpha-1}&\frac{z^2}{t^2}\frac{H_\nu^{(1)}(z)}{J_\nu(z)}\sqrt{zx}J_\nu(zx)\sqrt{zy}J_\nu(zy) dz\\
    &=c_1\int_{S_{\varepsilon,t}^+} \left(1-\frac{z^2}{t^2}\right)^{\alpha-1}\frac{z^2}{t^2}\frac{H_\nu^{(1)}(z)}{J_\nu(z)}\cos((x+y)z-2\theta_\nu)dz\\
    &\quad +c_2\int_{S_{\varepsilon,t}^+} \left(1-\frac{z^2}{t^2}\right)^{\alpha-1}\frac{z^2}{t^2}\frac{H_\nu^{(1)}(z)}{J_\nu(z)}\cos((x-y)z)dz\\
    &\quad +c_3\frac{x+y}{xy}\int_{S_{\varepsilon,t}^+}  \left(1-\frac{z^2}{t^2}\right)^{\alpha-1}\frac{z^2}{t^2}\frac{H_\nu^{(1)}(z)}{J_\nu(z)}\frac{\sin((x+y)z-2\theta_\nu)}{z}dz\\
    &\quad +c_4\frac{x-y}{xy}\int_{S_{\varepsilon,t}^+}  \left(1-\frac{z^2}{t^2}\right)^{\alpha-1}\frac{z^2}{t^2}\frac{H_\nu^{(1)}(z)}{J_\nu(z)}\frac{\sin((x-y)z)}{z}dz\\
    &\quad +\int_{S_{\varepsilon,t}^+} \left(1-\frac{z^2}{t^2}\right)^{\alpha-1}\frac{z^2}{t^2}\frac{H_\nu^{(1)}(z)}{J_\nu(z)}K(x,y,z) dz\\
    &=\sum_{j=1}^{5}T_{\varepsilon, t;j}^\alpha (x,y).
\end{align*}

Assume now that $0<\frac{x}{2}\leq y\leq 2x$. Let us estimate $\|T_{\varepsilon,\,\cdot\, ;j}^\alpha(x,y)\|_{L^2\left(\left(\frac{1}{|x-y|},\infty\right),\frac{dt}{t}\right)}$, for $j=1,...,5$. 
We will use that $\left|\frac{H_\nu^{(1)}(z)}{J_\nu (z)}\right|\leq Ce^{-2{\rm Im}\;z}$, $z\in S_{\varepsilon, t}^+$, $t>0$ (see \cite[(5.11.4) and (5.11.6)]{Le}), and that, if $z=t+iu$, $t,u>0$, 
\begin{equation}\label{product}
\left|1-\frac{z^2}{t^2}\right|^{\alpha-1}\frac{|z|^2}{t^2}\leq \left(\frac ut+\frac{u^2}{t^2}\right)^{\alpha-1} \left(1+\frac{u}{t}\right)^2=\left(\frac ut\right)^{\alpha -1}\left(1+\frac{u}{t}\right)^{\alpha +1}\leq \left(\frac ut\right)^{\alpha -1}+\left(\frac ut\right)^{2\alpha}. 
\end{equation}

We first study the cases $j=3,5$. Let us observe that, since $x\geq |x-y|$, then, if $t|x-y|\geq 1$,
$$
|K(x,y,z)|\bigg|_{z=t+iu}\leq C\frac{e^{(x+y)u}}{x^2t^2}\left(1+\frac{1}{xt}+\frac{1}{x^2t^2}\right)\leq C\frac{e^{(x+y)u}}{x^2t^2}\leq C\frac{e^{(x+y)u}}{xt},\quad u>0.
$$
By taking into account that $2-r-s>|r-s|$, $r,s\in (0,1)$, we get, for every $t\geq \frac{1}{|x-y|}$, 
\begin{align*}
   |T_{\varepsilon,t; 3}^\alpha(x,y)|+ |T_{\varepsilon,t; 5}^\alpha(x,y)|&\leq 
    \frac{C}{xt}\int_0^\infty \left(\left(\frac ut\right)^{\alpha -1}+\left(\frac ut\right)^{2\alpha}\right)e^{-(2-x-y)u}du\leq \frac{C}{x}\int_0^\infty (v^{\alpha -1}+v^{2\alpha })e^{-t|x-y|v}dv\\
    &\leq \frac{C}{x}\left(\frac{1}{(t|x-y|)^\alpha}+\frac{1}{(t|x-y|)^{2\alpha+1}}\right)\leq \frac{C}{x(t|x-y|)^\alpha}.
    \end{align*}
 It follows that 
\begin{align*}
\int_{\frac{1}{|x-y|}}^\infty \left(|T_{\varepsilon,t; 3}^\alpha(x,y)|^2+|T_{\varepsilon,t;5}^\alpha (x,y)|^2\right)\frac{dt}{t}&\leq \frac{C}{x^2|x-y|^{2\alpha }}\int_{\frac{1}{|x-y|}}^\infty \frac{1}{t^{2\alpha +1}}\leq \frac{C}{x^2}.
\end{align*}

On the other hand, for $j=2,4$, by using that $2-|x-y|\geq x\geq |x-y|$, we can write, for every $t\geq \frac{1}{|x-y|}$,
\begin{align*}
    |T_{\varepsilon,t; 2}^\alpha(x,y)|+|T_{\varepsilon,t; 4}^\alpha(x,y)|&\leq 
    C\left(1+\frac{|x-y|}{x^2t}\right)\int_0^\infty \left(\left(\frac ut\right)^{\alpha -1}+\left(\frac ut\right)^{2\alpha}\right)e^{-(2-|x-y|)u}du\\
    &\leq Ct\int_0^\infty (v^{\alpha -1}+v^{2\alpha })e^{-txv}dv\leq C\frac{t}{(tx)^\alpha},
    \end{align*}
and, since $\alpha >1$, 
\begin{align*}
\int_{\frac{1}{|x-y|}}^\infty \left(|T_{\varepsilon,t;2}^\alpha (x,y)|^2+T_{\varepsilon,t;2}^\alpha (x,y)|^2\right)\frac{dt}{t}&\leq \frac{C}{x^{2\alpha}}\int_{\frac{1}{|x-y|}}^\infty \frac{1}{t^{2\alpha -1}}dt\leq C\frac{|x-y|^{2\alpha -2}}{x^{2\alpha }}\leq \frac{C}{x^2}.
\end{align*}

In the following, we handle the case $j=1$. By partial integration we get
\begin{align*}
    T_{\varepsilon, t;1}^\alpha(x,y)&= i\int_\varepsilon^\infty \left(1-\frac{z^2}{t^2}\right)^{\alpha-1}\frac{z^2}{t^2}\frac{H_\nu^{(1)}(z)}{J_\nu(z)}\cos((x+y)z-2\theta_\nu)\bigg|_{z=t+iu}du\\
    &=\int_\epsilon^\infty \left(1-\frac{z^2}{t^2}\right)^{\alpha-1}\frac{z^2}{t^2}\frac{H_\nu^{(1)}(z)}{J_\nu(z)}\frac{d}{du}\left[\frac{\sin((x+y)z-2\theta _\nu)}{x+y}\right]\bigg|_{z=t+iu}du\\
    &=\left.\left(1-\frac{z^2}{t^2}\right)^{\alpha-1}\frac{z^2}{t^2}\frac{H_\nu^{(1)}(z)}{J_\nu(z)}\frac{\sin((x+y)z-2\theta_\nu)}{x+y}\bigg\rvert_{z=t+iu}\right]_{u=\varepsilon}^{u\rightarrow +\infty}\\
    &\quad -\int_\varepsilon^\infty\frac{d}{du}\left[\left(1-\frac{z^2}{t^2}\right)^{\alpha-1}\frac{z^2}{t^2}\frac{H_\nu^{(1)}(z)}{J_\nu(z)}\bigg|_{z=t+iu}\right]\frac{\sin((x+y)z-2\theta_\nu)}{x+y}\bigg|_{z=t+iu}du,\quad t>0.
\end{align*}

We have that, for every $t>\frac12 $,
\begin{align}\label{integral1}
    \left|\left(1-\frac{z^2}{t^2}\right)^{\alpha-1}\frac{z^2}{t^2}\frac{H_\nu^{(1)}(z)}{J_\nu(z)}\frac{\sin((x+y)z-2\theta _\nu)}{x+y}\bigg\rvert_{z=t+iu}\right|&\leq \frac{C}{x} \left(\left(\frac{u}{t}\right)^{\alpha -1}+\left(\frac{u}{t}\right)^{2\alpha}\right)e^{-(2-x-y)u},\quad u>0.
\end{align}
Observe that the right-hand side of the inequality and, consequently, the left one, converge to 0 when $u\rightarrow+\infty$. Also, since $\left(H_\nu^{(1)}\right)'(z)J_\nu(z)-J_\nu'(z)H_\nu^{(1)}(z)=\frac{2i}{\pi z}$, for $z\neq 0$ (see \cite[p. 94]{CR1}), we obtain
\begin{align*}
    \frac{d}{du}\left[\left(1-\frac{z^2}{t^2}\right)^{\alpha-1}\frac{z^2}{t^2}\frac{H_\nu^{(1)}(z)}{J_\nu(z)}\bigg\rvert_{z=t+iu}\right]&\\&\hspace{-3cm}=i\left[-(\alpha-1)\left(1-\frac{z^2}{t^2}\right)^{\alpha-2}\frac{2z^3}{t^4}\frac{H_\nu^{(1)}(z)}{J_\nu(z)}\right.\\
    &\hspace{-3cm}\quad +\left.\left(1-\frac{z^2}{t^2}\right)^{\alpha-1}\frac{2z}{t^2}\frac{H_\nu^{(1)}(z)}{J_\nu(z)}+\frac{2i}{\pi}\left(1-\frac{z^2}{t^2}\right)^{\alpha-1}\frac{z}{t^2}\frac{1}{(J_\nu(z))^2}\right]\bigg\rvert_{z=t+iu}\\
    &\hspace{-3cm}=\sum_{j=1}^3 B_j(z,t)\bigg|_{z=t+iu}, \quad t, u>0.
\end{align*}
Then, we can write
\begin{align*}
    T_{\varepsilon, t;1}^\alpha(x,y)&=  -\left(1-\frac{z^2}{t^2}\right)^{\alpha-1}\frac{z^2}{t^2}\frac{H_\nu^{(1)}(z)}{J_\nu(z)}\frac{\sin((x+y)z-2\theta_\nu)}{x+y}\bigg\rvert_{z=t+i\varepsilon}\\
    &\quad -\int_\varepsilon^\infty \sum_{j=1}^3B_j(z,t)\frac{\sin((x+y)z-2\theta_\nu)}{x+y}\bigg\rvert_{z=t+iu}du,\quad t>0.
\end{align*}
From \eqref{integral1}, since $\alpha >1$ and $0<\varepsilon <1$, it follows that,
$$
\left|\left(1-\frac{z^2}{t^2}\right)^{\alpha-1}\frac{z^2}{t^2}\frac{H_\nu^{(1)}(z)}{J_\nu(z)}\frac{\sin((x+y)z-2\theta_\nu)}{x+y}\bigg\rvert_{z=t+i\varepsilon}\right|\leq \frac{C}{x}\left(\frac{1}{t^{\alpha -1}}+\frac{1}{t^{2\alpha}}\right)\leq \frac{C}{xt^{\alpha -1}},\quad t>\frac12 .
$$
On the other hand, by proceeding as above, taking into account again that $2-r-s>|r-s|$, $r,s\in (0,1)$, and that $\alpha >1$, we get
\begin{align*}
    \left|\int_\varepsilon ^\infty \sum_{j=1}^3B_j(z,t)\frac{\sin((x+y)z-2\theta_\nu)}{x+y}\bigg\rvert_{z=t+iu}du\right|&\\
    &\hspace{-5cm}\leq \frac{C}{xt}\int_0^\infty \left(\left(\frac{u}{t}\right)^{\alpha -2}\left(1+\frac{u}{t}\right)^{\alpha +1}+\left(\frac{u}{t}\right)^{\alpha -1}\left(1+\frac{u}{t}\right)^\alpha\right)e^{-u(2-x-y)}du\\
    &\hspace{-5cm}\leq \frac{C}{xt}\int_0^\infty \left(\left(\frac{u}{t}\right)^{\alpha -2}+\left(\frac{u}{t}\right)^{\alpha -1}+\left(\frac{u}{t}\right)^{2\alpha-1}\right)e^{-u(2-x-y)}du\\
    &\hspace{-5cm}\leq \frac{C}{x}\int_0^\infty (v^{\alpha -2}+v^{\alpha -1}+v^{2\alpha -1})e^{-t|x-y|v}dv\\
    &\hspace{-5cm}\leq \frac{C}{x}\left(\frac{1}{(t|x-y|)^{\alpha -1}}+\frac{1}{(t|x-y|)^\alpha }+\frac{1}{(t|x-y|)^{2\alpha }}\right)\\
    &\hspace{-5cm}\leq \frac{C}{x(t|x-y|)^{\alpha -1}},\quad t>\frac12 .
\end{align*}
Thus, 
$$
|T_{\varepsilon,t;1}^\alpha (x,y)|\leq \frac{C}{x}\left(\frac{1}{t^{\alpha -1}}+\frac{1}{(t|x-y|)^{\alpha -1}}\right)\leq \frac{C}{x(t|x-y|)^{\alpha -1}},\quad t>\frac12 ,
$$
and, since $t>\frac12 $ when $t|x-y|\geq 1$, we deduce that 
$$
\int_{\frac{1}{|x-y|}}^\infty |T_{\varepsilon,t;1}^\alpha (x,y)|^2\frac{dt}{t}\leq \frac{C}{x^2|x-y|^{2\alpha -2}}\int_{\frac{1}{|x-y|}}^\infty \frac{1}{t^{2\alpha -1}}dt\leq \frac{C}{x^2}.
$$
Putting together all the estimates above, we conclude that \eqref{claim} is satisfied when $\frac{x}{2}\leq y\leq 2x$.

Next, we show \eqref{claim} for $y\in (0,\frac{x}{2})\cup (2x,1)$. Actually, it is sufficient to see the estimation when $0<y<\frac{x}{2}$, because of the symmetry of the property. 

Suppose that $0<y<\frac{x}{2}$. We will use that $|J_\nu (w)|\leq C|w|^\nu$ for $|w|\leq 1$ and $|{\rm arg}\;(w)|\leq \pi$ (\cite[(5.4.3]{Le}), and that, for every $r\in (0,\pi)$ there exists $C>0$ such that $|\sqrt{w}J_\nu (w)|\leq  Ce^{{\rm Im}\;w}$ for $|w|\geq 1$ and $|{\rm arg} \;w|< \pi -r$ (\cite[p. 199]{Wat}). 

Observe that when $z=t+iu\in S_{\varepsilon,t}^+$ and $t|x-y|\geq 1$, then $|xz|\geq xt>t|x-y|\geq 1$ and, consequently, $|\sqrt{xz}J_\nu (xz)|\leq Ce^{xu}$. By using \eqref{product} we can write, for $t\geq \frac{1}{|x-y|}$,
\begin{align*}
   \left|\int_{S_{\varepsilon,t}^+} \left(1-\frac{z^2}{t^2}\right)^{\alpha-1}\frac{z^2}{t^2}\frac{H_\nu^{(1)}(z)}{J_\nu(z)}\sqrt{zx}J_\nu(zx)\sqrt{zy}J_\nu(zy) dz\right|&\\
   &\hspace{-7cm} \leq C\int_0 ^\infty \left[\left(\frac ut\right)^{\alpha -1}+\left(\frac ut\right)^{2\alpha}\right]e^{-u(2-x)}|\sqrt{yz}J_\nu (yz)|\bigg\rvert_{z=t+iu}du \\
   &\hspace{-7cm}=C\left(\int_{u>0,\, y|t+iu|\geq 1 }+\int_{u>0,\, y|t+iu|<1}\right)\left[\left(\frac ut\right)^{\alpha -1}+\left(\frac ut\right)^{2\alpha}\right]e^{-u(2-x)}|\sqrt{yz}J_\nu (yz)|\bigg\rvert_{z=t+iu}du\\
   &\hspace{-7cm}=I_1(x,y;t)+I_2(x,y;t).
\end{align*}
We get, in the same way as before, and by taking into account that $\frac{x}{2}<x-y<x$, that
$$
I_1(x,y;t)\leq C\int_0^\infty \left[\left(\frac ut\right)^{\alpha -1}+\left(\frac ut\right)^{2\alpha}\right]e^{-u(2-x-y)}du\leq C\frac{t}{(tx)^\alpha },\quad t\geq \frac{1}{|x-y|}.  
$$
On the other hand, when $\nu \geq -\frac{1}{2}$, we have
$$
I_2(x,y;t)\leq C\int_0^\infty \left[\left(\frac ut\right)^{\alpha -1}+\left(\frac ut\right)^{2\alpha}\right]e^{-u(2-x)}du\leq C\frac{t}{(t(2-x))^\alpha}\leq C\frac{t}{(tx)^\alpha} ,\quad t\geq \frac{1}{|x-y|},  
$$
and, in the case that $-1<\nu <-\frac{1}{2}$, since $\alpha >1$, it follows that, for every $t\geq\frac{1}{|x-y|}$,
\begin{align*}
I_2(x,y;t)&\leq C\int_0^\infty \left[\left(\frac ut\right)^{\alpha -1}+\left(\frac ut\right)^{2\alpha}\right]e^{-u(2-x)}(yu)^{\nu +\frac{1}{2}}du\leq Ct^{\nu +\frac{3}{2}}y^{\nu +\frac{1}{2}}\int_0^\infty (v^{\alpha -1}+v^{2\alpha})e^{-t(2-x)v}v^{\nu +\frac{1}{2}}dv\\
&\leq Ct^{\nu +\frac{3}{2}}y^{\nu +\frac{1}{2}}\left(\frac{1}{(t(2-x))^{\alpha+\nu +\frac{1}{2}}}+\frac{1}{(t(2-x))^{2\alpha+\nu +\frac{3}{2}}}\right)\leq C\frac{t^{\nu +\frac{3}{2}}y^{\nu +\frac{1}{2}}}{(tx)^{\alpha +\nu +\frac{1}{2}}}=C\frac{y^{\nu +\frac{1}{2}}}{t^{\alpha-1}x^{\alpha +\nu +\frac{1}{2}}}.  
\end{align*}

Then, when $\nu\geq -\frac{1}{2}$,
\begin{align*}
\int_{\frac{1}{|x-y|}}^\infty\left|\int_{S_{\varepsilon,t}^+}\left(1-\frac{z^2}{t^2}\right)^{\alpha -1}\frac{z^2}{t^2}\frac{H_\nu ^{(1)}(z)\sqrt{xz}J_\nu (xz)\sqrt{yz}J_\nu (yz)}{J_\nu (z)}dz\right|^2\frac{dt}{t}&\\
&\hspace{-2cm}\leq \frac{C}{x^{2\alpha}}\int_{\frac{1}{|x-y|}}^\infty \frac{1}{t^{2\alpha -1}}dt\leq C\frac{|x-y|^{2\alpha -2}}{x^{2\alpha}}\leq \frac{C}{x^2}.
\end{align*}
If $-1<\nu<-\frac{1}{2}$, we obtain
\begin{align*}
\int_{\frac{1}{|x-y|}}^\infty\left|\int_{S_{\varepsilon,t}^+}\left(1-\frac{z^2}{t^2}\right)^{\alpha -1}\frac{z^2}{t^2}\frac{H_\nu ^{(1)}(z)\sqrt{xz}J_\nu (xz)\sqrt{yz}J_\nu (yz)}{J_\nu (z)}dz\right|^2\frac{dt}{t}&\leq C\left(\frac{1}{x^2}+\frac{y^{2\nu +1}}{x^{2\alpha +2\nu +1}}\int_{\frac{1}{|x-y|}}^\infty \frac{1}{t^{2\alpha -1}}dt\right)\\
&\leq C\left(\frac{1}{x^2}+\frac{y^{2\nu +1}}{x^{2\nu +3}}\right) \leq C\frac{y^{2\nu +1}}{x^{2\nu +3}}.
\end{align*}
Thus we have established that \eqref{claim} holds when $0<y<\frac{x}{2}$ or $2x<y<1$.

We conclude that, for each $x\in (0,1)$, 
$$
G\mathbb{A}_{\nu;2}^{\alpha,+}(f)(x)\leq C\left\{
\begin{array}{ll}
\displaystyle \frac{1}{x}\int_0^x|f(y)|dy+\int_x^1\frac{|f(y)|}{y}dy,&\displaystyle \nu \geq -\frac12 ,\\
&\\
\displaystyle \frac{1}{x^{\nu +3/2}}\int_0^x|f(y)|y^{\nu +\frac{1}{2}}dy+x^{\nu +\frac{1}{2}}\int_x^1\frac{|f(y)|}{y^{\nu +3/2}}dy,
&\displaystyle -1<\nu < -\frac12 .
\end{array}
\right.
$$
From Hardy's inequality (\cite{Muck}) we deduce that  $G\mathbb{A}_{\nu;2}^{\alpha,+}$ is bounded on $L^p(0,1)$, when $\alpha >1$ and any of the conditions ($i'$) or ($ii'$) is satisfied. 

\noindent $\bullet$ {\bf $L^p$--boundedness properties for $G\mathbb{A}_{\nu ;3}^\alpha$}. We first observe that
$$
|A_{\nu , t;3}^\alpha (x,y)|=\phi (t|x-y|)\left| t\frac{\partial}{\partial t}R_t^{\alpha ,L_\nu }(x,y)-H_{\nu ,t;1}^\alpha (x,y)\right|, \quad x,y\in (0,1), t>0.
$$
According to \cite[pp. 106 and 107]{CR1} we obtain, for every $x,y\in (0,1)$ and $t>0$,
 \begin{align*}
 \left|t\frac{\partial}{\partial t}R_t^{\alpha ,L_\nu }(x,y)\right|=2\alpha |(R_t^{\alpha -1,L_\nu}(x,y)-R_t^{\alpha ,L_\nu }(x,y))|\leq C\left\{
 \begin{array}{ll}
 (xy)^{\nu +\frac{1}{2}}t^{2(\nu +1)},&(x,y)\in Q_{t,1},\\
 &\\
 t,&(x,y)\in Q_{t,2},
 \end{array}
 \right.
 \end{align*}
 where, for each $t>0$, $$Q_{t,1}=\left\{(x,y): x,y\in (0,1) \mbox{ and }0<\max\{x,y\}\leq \frac{4}{t}\right\},$$ $$Q_{t,2}=\left\{(x,y): x,y\in (0,1) \mbox{ and }
 \frac{4}{t}\leq \max\{x,y\},\;|x-y|\leq \frac{2}{t}\right\}.$$
 Note that when $0<t\leq 4$, $Q_{t,1}=(0,1)\times (0,1)$ and $Q_{t,2}=\emptyset$, as illustrated below.

%\begin{figure}[h!]
%\begin{tikzpicture}[scale=1.6]
%    \draw[->] (-0.5,0) -- (2.5,0);
%    \draw[->] (0,-0.5) -- (0,2.5);
%    \draw (2.4,-0.25) node {$x$};
%    \draw (-0.25,2.4) node {$y$};
%    
%	\fill[color=gray]
%       (0,0) -- (0,1.3) -- (1.3,1.3) -- (1.3,0);
%    \fill[color=lightgray]
%	 	  --  plot [domain=1.3:1] (\x,\x-0.65)
%	 		(2,1.35) -- (2,2) -- (1.35,2)-- (1.3,1.95) -- (1.3,0.65);
%	\fill[color=lightgray]
%		 --  plot [domain=0.65:1.3] (\x,\x+0.65)
%		 (1.35,2) -- (1.35,1.3) -- (0.65,1.3);
%		 
%    \draw[-] (0,2) -- (2,2);
%    \draw[-] (2,0) -- (2,2);
%	\draw (0.8,0.8) node {$Q_{t,1}$};
%	\draw (1.55,1.55) node {$Q_{t,2}$};
%	\draw (2.55,1.1) node {$y=\frac{x}{2}$};
%	\draw (1.2,2.35) node {$y=2x$};
%	
%	\draw[dashed]  [domain=0:1.1]plot(\x,2*\x);  
%    \draw [dashed] [domain=0:2.2] plot(\x,0.5*\x);
%       
%    \draw[dashed]  [domain=0:1.35]plot(\x,\x+0.65);  
%	\draw[dashed]  [domain=0.65:2]plot(\x,\x-0.65); 
%	
%	\draw[-] (1,-0.05) -- (1,0.05);
% 	\draw (1,-0.3) node {$\frac{1}{2}$};
% 	\draw[-] (2,-0.05) -- (2,0.05);
% 	\draw (2,-0.25) node {$1$};
%    \draw[-] (0.65,-0.05) -- (0.65, 0.05);
%    \draw (0.65,-0.3) node {$\frac{2}{t}$};
%    \draw[-] (1.3,-0.05) -- (1.3, 0.05);
%    \draw (1.3,-0.3) node {$\frac{4}{t}$};
%    
%    \draw[-] (-0.05,1) -- (0.05,1);
%    \draw (-0.25,1) node {$\frac{1}{2}$};
%    \draw[-] (-0.05,2) -- (0.05,2);
%    \draw (-0.25,2) node {$1$};
%    \draw[-] (-0.05,0.65) -- (0.05,0.65);
%    \draw (-0.25,0.65) node {$\frac{2}{t}$};
%    \draw[-] (-0.05,1.3) -- (0.05, 1.3);
%    \draw (-0.25,1.3) node {$\frac{4}{t}$};
%\end{tikzpicture}
%
%\end{figure}

\begin{figure}[!htb]
    \centering
                 \includegraphics{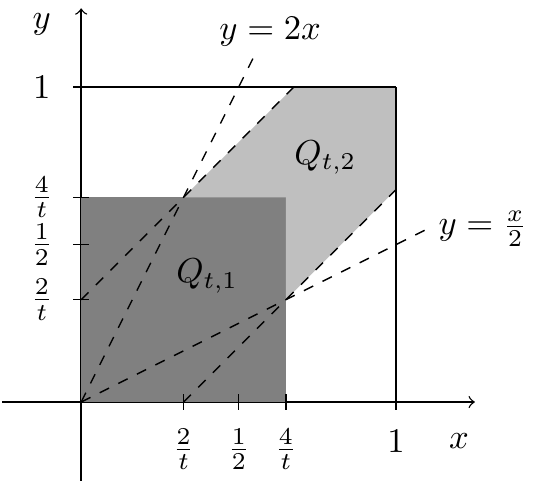}
              \end{figure}
 On the other hand, since $z^{-\nu }J_\nu (z)$ is bounded in $(0,1)$, we get, for every $t>0$ and $(x,y)\in Q_{t,1}$,
$$
 |H_{\nu ,t;1}^\alpha (x,y)|=2\alpha \sqrt{xy}t^2\left|\int_0^1(1-z^2)^{\alpha -1}z^3J_\nu(txz)J_\nu (tyz)dz\right|\leq C(xy)^{\nu +\frac{1}{2}}t^{2(\nu +1)}.
$$
 If $t>0$ and $(x,y)\in Q_{t,2}$, we take into account also that  $\sqrt{z}J_\nu (z)$ is bounded on $(1,\infty )$ and write
\begin{align*}
|H_{\nu ,t;1}^\alpha (x,y)|&\leq C\sqrt{xy}t^2\left(\int_0^{\frac{1}{tx}}+\int_{\frac{1}{tx}}^1\right) (1-z^2)^{\alpha -1}z^3|J_\nu (ztx)J_\nu (zty)|dz\\
& \leq C\sqrt{xy}t^2\left(\int_0^{\frac{1}{tx}}(1-z^2)^{\alpha-1}z^3(ztx)^\nu (zty)^\nu
dz+\int_{\frac{1}{tx}}^1(1-z^2)^{\alpha -1} z^3(txz)^{-\frac{1}{2}}(tyz)^{-\frac{1}{2}}dz\right) \\
&\leq C\left((xy)^{\nu +\frac{1}{2}}t^{2(\nu +1)}\int_0^{\frac{1}{tx}}(1-z^2)^{\alpha -1} z^{2\nu +3}dz+ t\int_0^1(1-z^2)^{\alpha -1}z^2dz\right)\\
& \leq Ct\left(\frac{y^{\nu +\frac{1}{2}}}{x^{\nu +\frac{1}{2}}}+1\right)\int_0^1(1-z^2)^{\alpha -1}z^2dz\leq Ct.
\end{align*}
In the last inequality we have used that, if $t>0$ and $(x,y)\in Q_{t,2}$ then $\frac{x}{2}<y<2x$.

Thus, since $\phi(t|x-y|)=0$, when $t|x-y|\geq 2$, we have seen that
 \begin{equation}\label{3.5}
 |A_{\nu , t;3}^\alpha (x,y)|\leq C\left\{
 \begin{array}{ll}
 (xy)^{\nu +\frac{1}{2}}t^{2(\nu +1)},& (x,y)\in  Q_{t,1},\\
 &\\
 t,& (x,y)\in (0,1)\times (0,1)\setminus Q_{t,1},
 \end{array} \right.
 \end{equation}
We note that, indeed, $A_{\nu,t;3}^\alpha (x,y)=0$, when $t>0$ and $(x,y)\in (0,1)\times (0,1)\setminus (Q_{t,1}\cup Q_{t,2})$ or $(x,y)\in Q_{t,1}\setminus\{(x,y):|x-y|\leq \frac{2}{t}\}$.
  
 On the other hand, since $\frac{d}{dx}[x^{-\nu }J_\nu (x)]=-x^{-\nu }J_{\nu +1}(x)$, $x\in (0,\infty )$, we can write, for each $ x,y\in (0,1)$ and $t>0$,
\begin{align*}
\frac{\partial}{\partial x}\left[t\frac{\partial}{\partial t}R_t^{\alpha ,L_\nu }(x,y)\right]&=\frac{\partial}{\partial x}\left(2\alpha\sum_{j=1}^\infty \left(1-\frac{s_{\nu ,j}^2}{t^2}\right)_+^{\alpha -1}\frac{s_{\nu ,j}^2}{t^2}\phi _{\nu ,j}(x)\phi _{\nu ,j}(y)\right)\\
&\hspace{-3cm}=2\alpha\sum_{j=1}^\infty \left(1-\frac{s_{\nu ,j}^2}{t^2}\right)_+^{\alpha -1}\frac{s_{\nu ,j}^2}{t^2}\phi _{\nu ,j}(y)\frac{\sqrt{2}}{\sqrt{s_{\nu ,j}}|J_{\nu +1}(s_{\nu ,j})|}\left(\frac{\nu +\frac12 }{x}\sqrt{s_{\nu ,j}x}J_\nu(s_{\nu ,j}x)-s_{\nu ,j}\sqrt{s_{\nu ,j}x}J_{\nu +1}(s_{\nu ,j}x)\right)\\
&\hspace{-3cm}=\frac{\nu +\frac{1}{2}}{x}\left[t\frac{\partial}{\partial t}R_t^{\alpha ,L_\nu }(x,y)\right]-2\alpha\sum_{j=1}^\infty \left(1-\frac{s_{\nu ,j}^2}{t^2}\right)_+^{\alpha -1}\frac{s_{\nu ,j}^3}{t^2}\phi _{\nu ,j}(y)\frac{\sqrt{2x}J_{\nu +1}(s_{\nu ,j}x)}{|J_{\nu +1}(s_{\nu ,j})|},
\end{align*}
and
\begin{align*}
\frac{\partial}{\partial x}H_{\nu ,t;1}^\alpha (x,y)&=2\alpha t\int_0^1(1-z)^{\alpha -1}z^2\left(\frac{\nu +\frac12 }{x}\sqrt{txz}J_\nu(txz)-tz\sqrt{txz}J_{\nu +1}(txz)\right)\sqrt{tyz}J_{\nu }(tyz)dz\\
&=\frac{(\nu +\frac{1}{2})t}{x}H_{\nu,t;1}^\alpha (x,y)-
2\alpha t^2\int_0^1(1-z)^{\alpha -1}z^3\sqrt{txz}J_{\nu +1}(txz)\sqrt{tyz}J_{\nu }(tyz)dz.
\end{align*}

By proceeding as above we obtain
\begin{equation}\label{3.8}
\left|\frac{\partial}{\partial x}\left[t\frac{\partial}{\partial t}R_t^{\alpha ,L_\nu }(x,y)\right]\right|+\left|\frac{\partial}{\partial x}H_{\nu ,t;1}^\alpha (x,y)\right|
\leq C\left\{
\begin{array}{ll}
\displaystyle \frac{(xy)^{\nu +\frac{1}{2}}t^{2(\nu +1)}}{x},&(x,y)\in Q_{t,1},\\
&\\
t^2,&(x,y)\in Q_{t,2},
\end{array}\right.
\end{equation}
with $x,y\in (0,1)$ and $t>0$.

Since for every $x,y\in (0,1)$ and $t>0$,  $|\frac{\partial}{\partial x}\phi (t|x-y|)|\leq Ct$, and $\phi (t|x-y|)=0$, when $t|x-y|\geq 2$, we deduce from \eqref{3.5} and \eqref{3.8} 
\begin{equation}\label{partial3.5}
 \left|\frac{\partial}{\partial x}A_{\nu , t;3}^\alpha (x,y)\right|\leq C\left\{
 \begin{array}{ll}
 \displaystyle \frac{(xy)^{\nu +\frac{1}{2}}t^{2(\nu +1)}}{x},& (x,y)\in  Q_{t,1},\\
 &\\
 t^2,& (x,y)\in (0,1)\times (0,1)\setminus Q_{t,1}.
 \end{array} \right.
 \end{equation}
 
We can see that
\begin{equation}\label{A3}
\|A_{\nu, \, \cdot\, ;3}^\alpha (x,y)\|_{\mathcal{H}}\leq \frac{C}{|x-y|}, \quad x,y\in (0,1),\; \frac{x}{2}<y<2x.
\end{equation}
Indeed, if $x,y\in (0,1)$, $\frac{x}{2}<y<2x$, then $2|x-y|<\max\{x,y\}$, and we can write, according to \eqref{3.5}, 
\begin{align*}
\Big\|A_{\nu, \, \cdot\, ;3}^\alpha (x,y)\Big\|_{\mathcal{H}}^2&=\left(\int_0^{\frac{4}{\max\{x,y\}}}+\int_{\frac{4}{\max\{x,y\}}}^\frac{2}{|x-y|}\right)|A_{\nu ,t; 3}^\alpha (x,y)|^2\frac{dt}{t}\leq C\left(x^{4\nu +2}\int_0^{\frac{4}{\max\{x,y\}}}t^{4\nu
+3}dt+\int_{\frac{4}{\max\{x,y\}}}^{\frac{2}{|x-y|}}tdt\right)\\
&\leq C\left(\frac{x^{4\nu +2}}{\max\{x,y\}^{4\nu +4}}+\frac{1}{|x-y|^2}+\frac{1}{\max\{x,y\}^2}\right)\leq C\left(\frac{1}{x^2}+\frac{1}{|x-y|^2}\right)\leq \frac{C}{|x-y|^2}.
\end{align*}

With the same arguments and using \eqref{partial3.5} we get
\begin{equation}\label{3.9}
\left\|\frac{\partial }{\partial x}A_{\nu, \, \cdot\, ;3}^\alpha (x,y)\right\|_{\mathcal{H}}+\left\|\frac{\partial }{\partial y}A_{\nu, \, \cdot\, ;3}^\alpha (x,y)\right\|_{\mathcal{H}}\leq \frac{C}{|x-y|^2},
\quad x,y\in (0,1), \;\frac{x}{2}<y<2x.
\end{equation}

On the other hand, when $0<y<\frac{x}{2}$ or $2x<y<1$ and $t|x-y|\leq 2$, we have that $(x,y)\in Q_{t,1}$, and, again from \eqref{3.5} we deduce
\begin{align*}
\Big\|A_{\nu, \, \cdot\, ;3}^\alpha (x,y)\Big\|_{\mathcal{H}}^2&=\int_0^\frac{2}{|x-y|}|A_{\nu ,t; 3}^\alpha (x,y)|^2\frac{dt}{t}\leq C(xy)^{2\nu +1}\int_0^\frac{2}{|x-y|}t^{4\nu
+3}dt=C\frac{(xy)^{2\nu +1}}{|x-y|^{4\nu +4}}\leq C\frac{(xy)^{2\nu +1}}{\max\{x,y\}^{4\nu +4}}.
\end{align*}
Thus,
\begin{equation}\label{A3global}
\Big\|A_{\nu, \, \cdot\, ;3}^\alpha (x,y)\Big\|_{\mathcal{H}}\leq C\frac{\min\{x,y\}^{\nu +\frac{1}{2}}}{\max\{x,y\}^{\nu +\frac{3}{2}}}\leq C\psi_\nu (x,y), \quad 0<y<\frac{x}{2} \mbox{ or }2x<y<1.
\end{equation}
Here, the function $\psi_\nu$ is the one appearing in \eqref{claim}. 

For every $x\in (0,1)$ we write
$$
G\mathbb{A}_{\nu ;3}^\alpha (f)(x)\leq G\mathbb{A}_{\nu ;3}^\alpha (\chi_{(\frac{x}{2},2x)}f)(x)+G\mathbb{A}_{\nu ;3}^\alpha (\chi_{(0,1)\setminus (\frac{x}{2},2x)}f)(x)=:G\mathbb{A}_{\nu ;3}^{\alpha,\rm loc} (f)(x)+G\mathbb{A}_{\nu ;3}^{\alpha,\rm{glob}} (f)(x).
$$

Suppose that $f\in L^2(0,1)$ and that the support of $f$ is compact in $(0,1)$. According to \eqref{3.5} we have that
$$
\int_0^1|A_{\nu ,t; 3}^\alpha (x,y)||f(y)|dy<\infty ,\quad x\in (0,1), \;t>0.
$$

Let $x\in (0,1)\setminus \supp (f)$ and $h\in \mathcal{H}$. We define
$$
F_x(t)=\left(\int_0^1A_{\nu, \, \cdot\, ; 3 }^\alpha (x,y)f(y)dy\right)(t),\quad t\in (0,\infty),
$$
where the integral is understood in an $\mathcal{H}$--Bochner sense. Note that by using \eqref{A3} and \eqref{A3global} we deduce
\begin{equation}\label{interchange}
\int_0^1\|A_{\nu, \, \cdot\, ;3}^\alpha (x,y)\|_{\mathcal{H}}|f(y)|dy<\infty.
\end{equation}

We can write 
$$
\int_0^\infty F_x(t)h(t)\frac{dt}{t}=\int_0^1\int_0^\infty A_{\nu ,t; 3}^\alpha(x,y)h(t)\frac{dt}{t}f(y)dy=\int_0^\infty \int_0^1A_{\nu ,t; 3}^\alpha (x,y)f(y)dyh(t)\frac{dt}{t}.
$$
The first equality follows from properties of the $\mathcal{H}$--Bochner integrals and the second one by interchanging the order of integration. This interchange can be justified by using H\"older inequality and \eqref{interchange}. 

Hence, 
$$
F_x(t)=\int_0^1A_{\nu ,t; 3}^\alpha (x,y)f(y)dy=\mathbb{A}_{\nu; 3} ^\alpha (f)(x,t),\quad \mbox{ a.e. } t\in (0,\infty ).
$$

According to local and vector--valued Calder\'on-Zygmund theory for singular integrals (see \cite{NS} for the scalar case) we can prove that the operator $\mathbb{A}_{\nu; 3} ^{\alpha, \rm loc}$ is bounded from $L^p(0,1)$ into $L^p((0,1),\mathcal{H})$, for every $1<p<\infty$,  and $\alpha >1$. In other words, we have seen that the operator $G\mathbb{A}_{\nu; 3} ^{\alpha ,{\rm loc}}$ is bounded on $L^p(0,1)$ for every $1<p<\infty$ and $\alpha >1$.

Moreover, by using the Minkowski integral inequality and estimation \eqref{A3global} we get, proceeding analogously to the analysis of $G\mathbb{A}_{\nu;2}^{\alpha, +}$, that $G\mathbb{A}_{\nu ;3}^{\alpha,\rm{glob}}$ is bounded on $L^p(0,\infty)$ when $\alpha >1$ and one of the conditions $(i')$ or $(ii')$ holds.

We conclude then that $G\mathbb{A}_{\nu ;3}^\alpha$ is a bounded operator on $L^p(0,\infty)$ when $\alpha >1$ and one of the conditions $(i')$ or $(ii')$ is satisfied.

From \eqref{decomposition} and the $L^p$--boundedness properties for $G\mathbb{A}_{\nu ;2}^\alpha$ and $G\mathbb{A}_{\nu ;3}^\alpha$ we obtain that the operator defined by
$$
G\mathbb{B}_\nu ^\alpha(f)(x)=\left\|\int_0^1H_{\nu, \cdot, ; 2}^\alpha (x,y)f(y)dy\right\|_{\mathcal{H}},\quad x\in (0,1),
$$
is bounded on $L^p(0,1)$, when $\alpha >1$, and provided that any of the conditions $(i')$ or $(ii')$ is fulfilled.

Now, we observe that, since $\{\phi _{\nu ,j}\}_{j=1}^\infty$ is an orthonormal complete sequence in $L^2(0,1)$, it follows that, when $\alpha >\frac{1}{2}$,
\begin{align*}
\int_0^1|G_\alpha ^{L_\nu} (f)(x)|^2dx&=4\alpha ^2\int_0^1\int_0^\infty \left|\sum_{j=1}^\infty \left(1-\frac{s_{\nu,j}^2}{t^2}\right)_+^{\alpha -1}\left(\frac{s_{\nu ,j}}{t}\right)^2c_{\nu, j}(f)\phi _{\nu ,j}(x)\right|^2\frac{dt}{t}dx\\
&=4\alpha ^2\int_0^\infty \sum_{j=1}^\infty \left(1-\frac{s_{\nu ,j}^2}{t^2}\right)_+^{2\alpha -2}\left(\frac{s_{\nu,j}}{t}\right)^4|c_{\nu ,j}(f)|^2\frac{dt}{t}\\
&=4\alpha ^2\int_0^\infty \sum_{j=1}^\infty (1-u^2)_+^{2(\alpha -1)}u^3|c_{\nu ,j}(f)|^2du=2\alpha ^2\frac{\Gamma (2\alpha -1)}{2\Gamma (2\alpha +1)}\|f\|_2^2.
\end{align*}

We have also seen that $G\mathbb{A}_{\nu; 1}^\alpha$ is bounded on $L^2(0,1)$ when $\alpha >\frac{1}{2}$. Then, again from \eqref{decomposition}, we deduce that $G\mathbb{B}_\nu ^\alpha$ is bounded on $L^2(0,1)$ provided that $\alpha >\frac{1}{2}$.

By using now Stein's complex interpolation theorem (see, for instance, \cite{St2}), we obtain that $G\mathbb{B}_\nu^\alpha$ is bounded on $L^p(0,1)$ when some of the properties $(i)$, $(ii)$ or $iii)$ in the statement of the Theorem \ref{Theorem1.1} holds. This fact, jointly with the $L^p$--boundedness properties of the operator $G\mathbb{A}_{\nu ;1}^\alpha $, leads to \eqref{Objective1}. 

Next, we are going to see that, for every $f\in L^p(0,1)$,
\begin{equation}\label{3.10}
\|f\|_p\leq C\|G_\alpha ^{L_\nu }(f)\|_p,
\end{equation}
provided that some of the previous properties $(i)$, $(ii)$ or $(iii)$ holds.

Note firstly that if $1<p<\infty$ and $p'=p/(p-1)$ is the conjugate of $p$ we have that
\begin{equation}\label{ii-p'}
\frac{2}{2\nu +3}<p\leq 2 \iff 2\leq p'<-\frac{2}{2\nu +1},    
\end{equation}
and
\begin{equation}\label{ii-p'-2}
    \alpha >\frac{1}{4\nu +4}\left(\frac{2}{p}+2\nu +1\right)\iff \alpha >\frac{1}{4\nu +4}\left(-\frac{2}{p'}+2\nu +3\right),
\end{equation}
when $-1<\nu <-\frac{1}{2}$. These relations allow us to prove \eqref{3.10} by using the $L^p$--boundedness properties for $G_\alpha ^{L_\nu}$ established in the first part of this proof and duality arguments.

Assume that some of the properties $(i)$, $(ii)$ or $(iii)$ is satisfied and suppose that $f\in C_c^\infty (0,1)$, the space of the smooth functions with compact support in $(0,1)$. We write 
$$
\mathbb{S}_{\nu ,t}=t\partial _tR_t^{\alpha ,L_\nu },\quad t>0.
$$

We have that
\begin{equation}\label{3.11}
f=c_\alpha \int_0^\infty (\mathbb{S}_{\nu ,t})^2f\frac{dt}{t},\quad \mbox{ with }c_\alpha =2-\frac{1}{\alpha},
\end{equation}
where the integral is understood in the $L^2(0,1)$--Bochner sense. Indeed, we can write
$$
(\mathbb{S}_{\nu ,t})^2f=4\alpha^2\sum_{j=1}^\infty \left(1-\frac{s_{\nu,j}^2}{t^2}\right)_+^{2(\alpha -1)}\frac{s_{\nu ,j}^4}{t^4}c_{\nu,j}(f)\phi _{\nu ,j},\quad t>0.
$$
 Partial integration leads to
\[c_{\nu,j}(f)=s_{\nu,j}^{-2} c_{\nu,j}(L_\nu f), \quad j\in \mathbb N,\]
which gives that $|c_{\nu,j}(f)|=O( s_{\nu,j}^{-2})$.

For every $n\in \mathbb{N}$, we have
\begin{align*}
\left\|\sum_{j=1}^\infty \left(1-\frac{s_{\nu,j}^2}{t^2}\right)_+^{2(\alpha -1)}\frac{s_{\nu ,j}^4}{t^4}c_{\nu,j}(f)\phi _{\nu ,j}\right\|_2&\leq \sum_{j=1}^\infty \left(1-\frac{s_{\nu,j}^2}{t^2}\right)_+^{2(\alpha -1)}\frac{s_{\nu ,j}^4}{t^4}|c_{\nu,j}(f)|\|\phi _{\nu ,j}\|_2\\
&\leq \sum_{j=1}^\infty \left(1-\frac{s_{\nu,j}^2}{t^2}\right)_+^{2(\alpha -1)}\frac{s_{\nu ,j}^4}{t^4}\frac{1}{s_{\nu,j}^2},\quad t>0.
\end{align*}

Since $s_{\nu ,j}\sim j$, for $j\in \mathbb{N}$ (\cite[p. 489]{Wat}), and $\alpha >\frac{1}{2}$, we obtain
$$
\int_0^\infty \sum_{j=1}^\infty \left(1-\frac{s_{\nu,j}^2}{t^2}\right)_+^{2(\alpha -1)}\frac{s_{\nu ,j}^4}{t^4}\frac{1}{s_{\nu ,j}^2}\frac{dt}{t}\leq C \sum_{j=1}^\infty \frac{1}{j^2}\int_0^1 (1-u^2)^{2(\alpha -1)}u^3du<\infty. 
$$

According to the Dominated Convergence Theorem for $L^2(0,1)$--Bochner integrals we can write
\begin{align*}
\int_0^\infty (\mathbb{S}_{\nu ,t})^2f\frac{dt}{t}&=4\alpha ^2\sum_{j=1}^\infty c_{\nu ,j}(f)\phi _{\nu ,j}\int_0^\infty \left(1-\frac{s_{\nu,j}^2}{t^2}\right)_+^{2(\alpha -1)}\frac{s_{\nu ,j}^4}{t^4}\frac{dt}{t}\\
&=4\alpha ^2\sum_{j=1}^\infty c_{\nu ,j}(f)\phi _{\nu ,j}\int_0^1 (1-u^2)^{2(\alpha -1)}u^3du=\frac{\alpha}{2\alpha -1}f.
\end{align*}
Thus \eqref{3.11} is proved.

Let $g\in C_c^\infty (0,1)$. We have that
\begin{align*}
\int_0^1g(x)(\mathbb{S}_{\nu ,t})^2f(x)dx&=4\alpha ^2\sum_{j=1}^\infty c_{\nu ,j}(g)c_{\nu ,j}(f)\left(1-\frac{s_{\nu,j}^2}{t^2}\right)_+^{2(\alpha -1)}\frac{s_{\nu ,j}^4}{t^4}\\
&=4\alpha ^2\sum_{j=1}^\infty \left(1-\frac{s_{\nu,j}^2}{t^2}\right)_+^{\alpha -1}\frac{s_{\nu ,j}^2}{t^2}c_{\nu ,j}(g)\left(1-\frac{s_{\nu,j}^2}{t^2}\right)_+^{\alpha -1}\frac{s_{\nu ,j}^2}{t^2}c_{\nu ,j}(f)\\
&=\int_0^1\mathbb{S}_{\nu ,t}f(x)\mathbb{S}_{\nu ,t}g(x)dx.
\end{align*}
Then, the properties of the $L^2(0,1)$--Bochner integrals lead to
\begin{align*}
\int_0^1f(x)g(x)dx&=c_\alpha \int_0^1g(x)\left(\int_0^\infty (\mathbb{S}_{\nu ,t})^2f\frac{dt}{t}\right)(x)dx=c_\alpha \int_0^\infty \int_0^1g(x)(\mathbb{S}_{\nu ,t})^2f(x)dx\frac{dt}{t}\\
&=c_\alpha \int_0^\infty \int_0^1(\mathbb{S}_{\nu ,t}f)(x)(\mathbb{S}_{\nu ,t}g)(x)dx\frac{dt}{t}= c_\alpha \int_0^1\int_0^\infty(\mathbb{S}_{\nu ,t}f)(x)(\mathbb{S}_{\nu ,t}g)(x)\frac{dt}{t}dx.
\end{align*}
The last equality is justified because
$$
\int_0^1\int_0^\infty \left|(\mathbb{S}_{\nu ,t}f)(x)(\mathbb{S}_{\nu ,t}g)(x)\right|\frac{dt}{t}dx\leq \|G_\alpha ^{L_\nu }f\|_p\|G_\alpha ^{L_\nu }g\|_{p'}\leq C\|f\|_p\|g\|_{p'},
$$
where we have used that $G_\alpha ^{L_\nu}$ is bounded on $L^p(0,1)$ and on $L^{p'}(0,1)$ when some of the properties $(i)$, $(ii)$ or $(iii)$ holds (recall the relations \eqref{ii-p'} and \eqref{ii-p'-2}). 

Since $C_c^\infty (0,1)$ is a dense subspace of $L^p(0,1)$ we conclude that
$$
\int_0^1f(x)g(x)dx=c_\alpha \int_0^1\int_0^\infty (\mathbb{S}_{\nu ,t}f)(x)(\mathbb{S}_{\nu ,t} g)(x)\frac{dt}{t}dx,
$$
for every $f\in L^p(0,1)$ and $g\in L^{p'}(0,1)$.

Let $f\in L^p(0,1)$. It follows that
\begin{align*}
\|f\|_p&=\sup_{g\in L^{p'}(0,1), \|g\|_{p'}\leq 1}\left|\int_0^1f(x)g(x)dx\right|=c_\alpha \sup_{g\in L^{p'}(0,1), \|g\|_{p'}\leq 1}\left|\int_0^1\int_0^\infty (\mathbb{S}_{\nu ,t}f)(x)(\mathbb{S}_{\nu ,t} g)(x)\frac{dt}{t}dx\right|\\
&\leq c_\alpha \sup_{g\in L^{p'}(0,1), \|g\|_{p'}\leq 1}\|G_\alpha ^{L_\nu }f\|_p\|G_\alpha ^{L_\nu }g\|_{p'}\leq C\|G_\alpha ^{L_\nu }f\|_p,
\end{align*}
and \eqref{3.10} is established.

Thus, the proof of Theorem \ref{Theorem1.1} is finished.

\section{Proof of Theorem \ref{Theorem1.3}}

We adapt the ideas developed in \cite[pp. 118--120]{BS} to a Hilbert--valued setting. 

Suppose that $G_\alpha^{L_\nu}$ is bounded on $L^p(0,1)$. We are going to show that there exists $C>0$ such that
$$\|G_\alpha^{B_\nu}f\|_{L^p(0,\infty)}\leq C\|f\|_{L^p(0,\infty)},\;\;\;\;f\in C_c^\infty(0,\infty).$$

Let $f\in C_c^\infty(0,\infty)$, the space of smooth functions with compact support in $(0,\infty)$. Without loss of generality, we will assume that $\mbox{supp}(f)\subset (0,a)$, with $a>0$. For every $r>0$, we define $f_r(x)=f(rx)$, $x\in (0,\infty)$. Note that there exists $r_0>1$ such that, for every $r>r_0$, $f_r\in C_c^\infty(0,1)$, where $C_c^\infty(0,1)$ represents the space of smooth functions with compact support in $(0,1)$.

We recall that for $g\in C_c^\infty(0,1)$ and $h\in C_c^\infty(0,\infty)$
$$G_\alpha^{L_\nu}(g)(x)=2\alpha\left\|\int_0^1g(y)\sum_{j=1}^\infty\left(1-\frac{s_{\nu,j}^2}{t^2}\right)_+^{\alpha -1}\frac{s_{\nu,j}^2}{t^2}\phi_{\nu,j}(x)\phi_{\nu,j}(y)dy\right\|_{\mathcal H},\;\;\;x\in (0,1),$$

$$G_\alpha^{B_\nu}(h)(x)=2\alpha\left\|\int_0^\infty h(y)\int_0^\infty \left(1-\frac{z^2}{t^2}\right)_+^{\alpha -1}\frac{z^2}{t^2}\sqrt{xz}J_{\nu}(xz)\sqrt{yz}J_{\nu}(yz)dzdy\right\|_{\mathcal H},\;\;\;x\in (0,\infty).$$

Let $r>r_0$. We define 
\[F_r(x)=\chi_{(0,r)}(x) \sum_{j=1}^\infty \left(1-\frac{s_{\nu,j}^2}{t^2 r^2}\right)_+^{\alpha-1} \frac{s_{\nu,j}^2}{t^2 r^2} c_{\nu,j}(f_r)\phi_{\nu, j}\left(\frac xr\right).\]
In view of the boundedness of $G_\alpha^{L_\nu}$ on $L^p(0,1)$, by considering $\mathcal H=L^2\left((0,\infty),\frac{dt}{t}\right)$, we obtain 
\begin{align*}
    \|F_r\|_{L^p((0,\infty),\mathcal H)}&=\left\|\left(\int_0^\infty \left|\sum_{j=1}^\infty \left(1-\frac{s_{\nu,j}^2}{t^2 r^2}\right)_+^{\alpha-1} \frac{s_{\nu,j}^2}{t^2 r^2} c_{\nu,j}(f_r)\phi_{\nu, j}\left(\frac{\cdot}{r}\right)\right|^2\frac{dt}{t}\right)^{\frac{1}{2}}\right\|_{L^p(0,r)}\\
    &=r^{\frac{1}{p}}\left\|\left(\int_0^\infty \left|\sum_{j=1}^\infty \left(1-\frac{s_{\nu,j}^2}{t^2}\right)_+^{\alpha-1} \frac{s_{\nu,j}^2}{t^2} c_{\nu,j}(f_r)\phi_{\nu, j}\left(\cdot\right)\right|^2\frac{dt}{t}\right)^{\frac{1}{2}}\right\|_{L^p(0,1)}\\
    &\leq C r^{\frac{1}{p}}\|f_r\|_{L^p(0,1)}=C\|f\|_{L^p(0,\infty)}.
\end{align*}

Since Plancherel's equality holds in $L^2((0,1),\mathcal H)$ with respect to $\{\phi_{\nu,j}\}_{j=1}^\infty$, we also have that
\begin{align}\label{4.1}
    \|F_r\|_{L^2((0,\infty),\mathcal H)}&=r^{\frac{1}{2}}\left\|\left(\int_0^\infty \left|\sum_{j=1}^\infty \left(1-\frac{s_{\nu,j}^2}{t^2}\right)_+^{\alpha-1} \frac{s_{\nu,j}^2}{t^2} c_{\nu,j}(f_r)\phi_{\nu, j}\left(\cdot\right)\right|^2\frac{dt}{t}\right)^{\frac{1}{2}}\right\|_{L^2(0,1)}\\
    \nonumber&=r^{\frac{1}{2}}\left(\sum_{j=1}^\infty |c_{\nu,j}(f_r)|^2 \int_0^\infty  \left(1-\frac{s_{\nu,j}^2}{t^2}\right)_+^{2(\alpha-1)}\frac{s_{\nu,j}^4}{t^4}\frac{dt}{t}\right)^{\frac{1}{2}}\\
    \nonumber&=r^{\frac{1}{2}}\left(\sum_{j=1}^\infty |c_{\nu,j}(f_r)|^2 \int_0^1  \left(1-z^2\right)^{2(\alpha-1)}z^3 dz\right)^{\frac{1}{2}}\\
    \nonumber&=\left(\frac{\Gamma(2\alpha-1)}{2\Gamma(2\alpha+1)}\right)^{\frac{1}{2}} r^{\frac{1}{2}} \left(\sum_{j=1}^\infty |c_{\nu,j}(f_r)|^2 \right)^{\frac{1}{2}}\\
    \nonumber&=\left(\frac{\Gamma(2\alpha-1)}{2\Gamma(2\alpha+1)}\right)^{\frac{1}{2}}\|f\|_{L^2(0,\infty)}.
\end{align}
Then, there exist an increasing sequence $\{r_j\}_{j=1}^\infty\subset (r_0,\infty)$ such that $r_j\rightarrow \infty$ as $j\rightarrow \infty$, and a function $F\in L^p((0,\infty),\mathcal H)\cap L^2((0,\infty),\mathcal H)$ such that $F_{r_j}\rightarrow F$ as $j\rightarrow \infty$, weakly in $L^2((0,\infty),\mathcal H)$ and in $L^p((0,\infty),\mathcal H)$. Moreover,
\[\|F\|_{L^2((0,\infty),\mathcal H)}\leq C\|f\|_{L^2(0,\infty)},\]
and
\[\|F\|_{L^p((0,\infty),\mathcal H)}\leq C\|f\|_{L^p(0,\infty)}.\]

We are going to prove that
\[F(x,t)=\frac{1}{2\alpha} t\partial_t R_t^{\alpha, B_\nu}(f)(x), \quad \textrm{a.e. }(x,t)\in (0,\infty)^2.\]
Given $r>r_0$, for every $k\in \mathbb N$ we decompose
\[F_r=I_r^k+D_r^k,\]
where
\[I_r^k(x,t)=\chi_{(0,r)}(x)\sum_{j=1}^{k [r]} \left(1-\frac{s_{\nu,j}^2}{t^2 r^2}\right)_+^{\alpha-1} \frac{s_{\nu,j}^2}{t^2 r^2} c_{\nu,j}(f_r)\phi_{\nu, j}\left(\frac xr\right), \quad x,t\in (0,\infty),\]
and
\[D_r^k(x,t)=\chi_{(0,r)}(x)\sum_{j=k [r]+1}^\infty \left(1-\frac{s_{\nu,j}^2}{t^2 r^2}\right)_+^{\alpha-1} \frac{s_{\nu,j}^2}{t^2 r^2} c_{\nu,j}(f_r)\phi_{\nu, j}\left(\frac xr\right), \quad x,t\in (0,\infty).\]
We recall that 
\[c_{\nu,j}(f_r)=\left(\frac{r}{s_{\nu,j}}\right)^2 c_{\nu,j}((L_\nu f)_r), \quad j\in \mathbb N.\]
Let $k\in \mathbb N$. By proceeding as in \eqref{4.1}, we get
\begin{align*}
    \|D_r^k\|_{L^2((0,\infty),\mathcal H)}^2&=\frac{\Gamma(2\alpha-1)}{2\Gamma(2\alpha+1)} r \sum_{j=k[r]+1}^\infty |c_{\nu,j}(f_r)|^2\\
    &=\frac{\Gamma(2\alpha-1)}{2\Gamma(2\alpha+1)} r \sum_{j=k[r]+1}^\infty \frac{r^4}{s_{\nu,j}^4} |c_{\nu,j}((L_\nu f)_r)|^2\\
    &\leq \frac{Cr}{k^4}\sum_{j=1}^\infty |c_{\nu,j}((L_\nu f)_r)|^2=\frac{C}{k^4}\|L_\nu f\|_2^2,
\end{align*}
where the constant $C$ does not depend on $r$ or $k$. We have used that $s_{\nu,j}=\pi \left(j+\frac{\nu}{2}-\frac14+O\left(\frac1j\right)\right)$, for $j\in \mathbb N$ (see \cite[(2.2)]{BS}).

By using a diagonal argument we can find an increasing function $\phi:\mathbb N \rightarrow \mathbb N$ such that, for every $k\in \mathbb N$, there exists $D^k\in L^2((0,\infty),\mathcal H)$ verifying $D_{r_{\phi(j)}}^k\rightarrow D^k$, as $j\rightarrow \infty$, weakly in $L^2((0,\infty),\mathcal H)$, with $\|D^k\|_{L^2((0,\infty),\mathcal H)}\leq C/k^2$. Then, there exists an increasing function $\psi:\mathbb N \rightarrow \mathbb N$ such that $\|D^{\psi(k)}(x,\cdot)\|_{\mathcal H}\rightarrow 0$, as $k\rightarrow \infty$, for a.e. $x\in (0,\infty)$, and $D^{\psi(k)}(x,t)\rightarrow 0$, as $k\rightarrow \infty$, for a.e. $(x,t)\in (0,\infty)^2$. We define, for every $k\in \mathbb N$, $F^k=F-D^{\psi(k)}$. Hence, we have
$$
\lim_{j\rightarrow \infty}I_{r_{\phi(j)}}^{\psi(k)}=F_{r_{\phi(j)}}-D_{r_{\phi(j)}}^{\psi(k)}\longrightarrow F-D^{\psi(k)}=F^k,\mbox{ weakly in }L^2((0,\infty),\mathcal H).
$$  
Also, $F^k(x,t)\rightarrow F(x,t)$, as $k\rightarrow \infty$, for a.e. $(x,t)\in (0,\infty)^2$, and $\|F^k(x,\cdot)-F(x,\cdot)\|_{\mathcal H}\rightarrow 0$, as $k\rightarrow \infty$, for a.e. $x\in (0,\infty)$.

For every $k\in \mathbb{N}$ and $x\in (0,\infty)$ we denote by $\mathfrak{F}_x:[0,\pi k]\rightarrow \mathcal{H}$ the function given by
\[[\mathfrak F_x(y)](t)=\left(1-\frac{y^2}{t^2}\right)_+^{\alpha-1} \frac{y^2}{t^2} h_{\nu}(f)(y)\sqrt{xy} J_\nu(xy), \quad y\in (0,\pi k) \textrm{ and }t\in (0,\infty),\quad y\in (0,\pi k],\;t>0,
\]
and
$$
[\mathfrak F_x(0)](t)=0,\quad t>0.
$$

Our next objective is to establish that, for each $k\in \mathbb{N}$ and $x\in (0,\infty)$,
$$
\lim_{r \rightarrow \infty}I_r ^k(x,\cdot )=\int_0^{\pi k}\mathfrak{F}_x(y)dy,\quad \mbox{ in }\mathcal{H}.
$$
Here, the integral is understood in the $\mathcal{H}$--Bochner sense.

Let $k\in \mathbb N$ and $x\in(0,\infty)$. We first observe that, $\mathfrak{F}_x$ is $\mathcal{H}$--Bochner integrable. Since $J_\mu(z)\sim \frac{1}{2^\mu \Gamma(\mu+1)}z^{\mu}$, as $z\rightarrow 0^+$, when $\mu>-1$, we have that
\begin{align}\label{4.3}
    \nonumber\left\|\left(1-\frac{y^2}{t^2}\right)_+^{\alpha-1} \frac{y^2}{t^2} h_{\nu}(f)(y)\sqrt{xy} J_\nu(xy)\right\|_{\mathcal H}&\leq C \left\|\left(1-\frac{y^2}{t^2}\right)_+^{\alpha-1} \frac{y^2}{t^2}\right\|_{\mathcal{H}}|h_{\nu}(f)(y)| |\sqrt{xy} J_\nu(xy)|\\
    &\hspace{-2cm}\leq C\int_0^a |\sqrt{yz} J_\nu(yz)| |f(z)| dz |\sqrt{xy} J_\nu(xy)|\leq Cy^{2\nu+1}, \quad 0<y<\pi k+1.
\end{align}
Then,
$$
\int_0^{\pi k}\|\mathfrak{F}_x(y)\|_{\mathcal{H}}dy\leq \int_0^{\pi k}y^{2\nu +1}dy\leq C.
$$

Furthermore, ${\mathfrak F}_x$ is continuous in $[0,\pi k]$. In order to see this, assume that $y\in [0,\pi k]$ and $\{y_n\}_{n=1}^\infty$ is a sequence in $[0,\pi k]$ verifying $y_n\rightarrow y$, as $n\rightarrow \infty$. We are going to show that ${\mathfrak F}_x(y_n)\rightarrow {\mathfrak F}_x(y)$, as $n\rightarrow \infty$, in $\mathcal H$. Suppose firstly that $y>0$. We have that
\begin{align*}
    [\mathfrak F_x(y)](t)-[\mathfrak F_x(y_n)](t)&=\left(1-\frac{y^2}{t^2}\right)_+^{\alpha-1} \frac{y^2}{t^2}\left[ h_{\nu}(f)(y)\sqrt{xy} J_\nu(xy)-h_{\nu}(f)(y_n)\sqrt{xy_n} J_\nu(xy_n)\right]\\
    &\quad + \left[\left(1-\frac{y^2}{t^2}\right)_+^{\alpha-1} \frac{y^2}{t^2}-\left(1-\frac{y_n^2}{t^2}\right)_+^{\alpha-1} \frac{y_n^2}{t^2}\right]h_{\nu}(f)(y_n)\sqrt{xy_n} J_\nu(xy_n)\\
    &=H_1(t,n)+H_2(t,n), \quad n\in \mathbb N \textrm{ and }t>0.
\end{align*}
From the continuity of $h_\nu(f)$ on $[0,\infty)$, we get
\begin{align*}
    \|H_1(\cdot,n)\|_{\mathcal H} & =\left|h_{\nu}(f)(y)\sqrt{xy} J_\nu(xy)-h_{\nu}(f)(y_n)\sqrt{xy_n} J_\nu(xy_n)\right|\left\|\left(1-\frac{y^2}{t^2}\right)_+^{\alpha-1} \frac{y^2}{t^2}\right\|_{\mathcal{H}}\\
    &=\left(\frac{\Gamma(2\alpha-1)}{2\Gamma(2\alpha+1)}\right)^{\frac{1}{2}} \left|h_{\nu}(f)(y)\sqrt{xy} J_\nu(xy)-h_{\nu}(f)(y_n)\sqrt{xy_n} J_\nu(xy_n)\right|\underset{n\rightarrow\infty}{\longrightarrow} 0.
\end{align*}

On the other hand, since $|h_{\nu}(f)(y_n)\sqrt{xy_n} J_\nu(xy_n)|\leq Cy_n^{2\nu+1}$, $n\in \mathbb N$, we have that
$$
\|H_2(\cdot, n)\|_{\mathcal{H}}\leq Cy_n^{2\nu+1}\left\|\left(1-\frac{y^2}{t^2}\right)_+^{\alpha-1} \frac{y^2}{t^2}-\left(1-\frac{y_n^2}{t^2}\right)_+^{\alpha-1} \frac{y_n^2}{t^2}\right\|_{\mathcal{H}},\quad n\in \mathbb{N}.
$$

Let us write $m_n=\min\{y_n,y\}$ and $M_n=\max\{y_n,y\}$, $n\in \mathbb{N}$. Then, we get
\begin{align*}
     \left\|\left(1-\frac{y^2}{t^2}\right)_+^{\alpha-1} \frac{y^2}{t^2}-\left(1-\frac{y_n^2}{t^2}\right)_+^{\alpha-1} \frac{y_n^2}{t^2}\right\|_{\mathcal{H}}^2 \\
   &\hspace{-5cm} = \int_{m_n}^{M_n}\left|\left(1-\frac{m_n^2}{t^2}\right)^{\alpha -1}\frac{m_n^2}{t^2}\right|^2\frac{dt}{t}+\int_{M_n}^\infty \left|\left(1-\frac{m_n^2}{t^2}\right)^{\alpha-1} \frac{m_n^2}{t^2}-\left(1-\frac{M_n^2}{t^2}\right)^{\alpha-1} \frac{M_n^2}{t^2}\right|^2\frac{dt}{t}\\
   &\hspace{-5cm}=\int_{\frac{m_n}{M_n}}^1(1-u^2)^{2(\alpha -1)}u^3du+\int_0^1\left|\left(1-\frac{m_n^2}{M_n^2}u^2\right)^{\alpha-1}\frac{m_n^2}{M_n^2} -(1-u^2)^{\alpha-1} \right|^2u^3du, \quad n\in \mathbb N.
\end{align*}
Since $\alpha >\frac{1}{2}$ and $\frac{m_n}{M_n}\longrightarrow 1$, as $n\rightarrow \infty$, it is clear that the first integral converges to 0, as $n\rightarrow \infty$. On the other hand, we choose $n_0\in \mathbb{N}$ such that $\frac{m_n}{M_n}\in \left(\frac{1}{2},1\right)$, $n\geq n_0$, and observe that, for every $n\geq n_0$,
\[\left(1-\frac{m_n^2}{M_n^2}u^2\right)^{\alpha-1}\frac{m_n^2}{M_n^2}\leq \left\{ \begin{array}{ll}
   \left(1-u^2\right)^{\alpha-1}, & \frac12<\alpha<1, \\
   \left(1-\frac12u^2\right)^{\alpha-1}, & \alpha\geq 1.
\end{array}\right.\]
Then, by using the Dominated Convergence Theorem, we obtain that 
\[\lim_{n\rightarrow \infty} \int_0^1 \left|\left(1-\frac{m_n^2}{M_n^2}u^2\right)^{\alpha-1} \frac{m_n^2}{M_n^2}-(1-u^2)^{\alpha-1}\right|^2u^3du=0.\]
Hence, $\|H_2(\cdot,n)\|_{\mathcal H}\rightarrow 0$, as $n\rightarrow \infty$.

If $y=0$, we have that 
\begin{align*}
\left\|\mathfrak{F}_x(0)-\mathfrak{F}_x(y_n)\right\|_{\mathcal{H}}&=\left\|\left(1-\frac{y_n^2}{t^2}\right)_+^{\alpha -1}\frac{y_n^2}{t^2}\right\|_{\mathcal{H}}|h_\nu (f)(y_n) \sqrt{xy_n}J_\nu (xy_n)|\\
&= \left(\frac{\Gamma (2\alpha -1)}{2\Gamma (2\alpha +1)}\right)^{\frac{1}{2}}|h_\nu (f)(y_n)\sqrt{xy_n}J_\nu (xy_n)|\longrightarrow 0,\mbox{ as }n\rightarrow \infty.
\end{align*}

Since ${\mathfrak F}_x$ is continuous in $[0,\pi k]$ and $\mathcal{H}$ is separable, the $\mathcal{H}$--Bochner integral of ${\mathfrak F}_x$ on $[0,\pi k]$ coincides with the $\mathcal{H}$-Riemann integral of ${\mathfrak F}_x$ on $[0,\pi k]$ (see, for instance \cite{KSSTZ}).

Now we write
\begin{align*}
    I_r^k(x,t)&=\chi_{(0,r)}(x)\sum_{j=1}^{k [r]} \left(1-\frac{s_{\nu,j}^2}{t^2 r^2}\right)_+^{\alpha-1} \frac{s_{\nu,j}^2}{t^2 r^2} c_{\nu,j}(f_r)\phi_{\nu, j}\left(\frac xr\right)\\
    &=\chi_{(0,r)}(x)\sum_{j=1}^{k [r]} \left(1-\frac{s_{\nu,j}^2}{t^2 r^2}\right)_+^{\alpha-1} \frac{s_{\nu,j}^2}{t^2 r^2} h_{\nu}(f)\left(\frac{s_{\nu,j}}{r}\right)\sqrt{\frac xr s_{\nu,j}} J_\nu\left(\frac xr s_{\nu,j}\right)\frac{d_{\nu,j}^2}{r}\\
    &=\chi_{(0,r)}(x)\sum_{j=1}^{k [r]} \left(1-\frac{s_{\nu,j}^2}{t^2 r^2}\right)_+^{\alpha-1} \frac{s_{\nu,j}^2}{t^2 r^2} h_{\nu}(f)\left(\frac{s_{\nu,j}}{r}\right)\sqrt{\frac xr s_{\nu,j}} J_\nu\left(\frac xr s_{\nu,j}\right)\frac{s_{\nu, j+1}-s_{\nu,j}}{r}\\
     &\quad + \chi_{(0,r)}(x)\sum_{j=1}^{k [r]} \left(1-\frac{s_{\nu,j}^2}{t^2 r^2}\right)_+^{\alpha-1} \frac{s_{\nu,j}^2}{t^2 r^2} h_{\nu}(f)\left(\frac{s_{\nu,j}}{r}\right)\sqrt{\frac xr s_{\nu,j}} J_\nu\left(\frac xr s_{\nu,j}\right)\frac{d_{\nu,j}^2-s_{\nu, j+1}+s_{\nu,j}}{r}\\
    &=I_{r;1}^k(x,t)+I_{r;2}^k(x,t),\quad t>0.
\end{align*}
Here, $d_{\nu,j}=\sqrt{2}\,\left|\sqrt{s_{\nu,j}}J_{\nu+1}(s_{\nu,j})\right|$, $j\in \mathbb N$.

Let us show that
$$
\lim_{r\rightarrow \infty}I_{r;2}^k(x,\cdot)=0,\quad \mbox{ in }\mathcal{H}.
$$

Since $J_\mu(z)\sim \frac{1}{2^\mu \Gamma(\mu+1)}z^{\mu}$, as $z\rightarrow 0^+$, for any $\mu>-1$, and $s_{\nu,j}=\pi \left(j+\frac{\nu}{2}-\frac14+O\left(\frac1j\right)\right)$, for $j\in \mathbb N$, we get
\[\left|\sqrt{\frac xr s_{\nu,j}} J_\nu\left(\frac xr s_{\nu,j}\right)\right|\leq C \left(\frac xr s_{\nu,j}\right)^{\nu+\frac{1}{2}}\leq C \left(\frac{xj}{r}\right)^{\nu+\frac{1}{2}}, \quad j\in \mathbb N,\, 1\leq j\leq k[r],\]
and, 
\begin{align*}
    \left|h_{\nu}(f)\left(\frac{s_{\nu,j}}{r}\right)\right|&\leq C\int_0^a |f(y)| \left|\sqrt{\frac yr s_{\nu,j}} J_\nu\left(\frac yr s_{\nu,j}\right)\right| dy\\
    &\leq C \left(\frac{j}{r}\right)^{\nu+\frac{1}{2}}\int_0^a |f(y)|y^{\nu+\frac{1}{2}}dy, \quad j\in \mathbb N,\, 1\leq j\leq k[r].
\end{align*}
Also, since $d_{\nu,j}=\sqrt{\pi} \left(1+O\left(\frac1j\right)\right)$, for $j\in \mathbb N$, it follows that $d_{\nu,j}^2-s_{\nu, j+1}+s_{\nu,j}=O\left(\frac1j\right)$, $j\in \mathbb N$.

We then deduce  that 
\begin{align*}
        \|I_{r;2}^k(x,\cdot)\|_{\mathcal{H}}& \leq C\sum_{j=1}^{k [r]} \left\|\left(1-\frac{s_{\nu,j}^2}{t^2 r^2}\right)_+^{\alpha-1} \frac{s_{\nu,j}^2}{t^2 r^2}\right\|_{\mathcal H} \left(\frac{j}{r}\right)^{2\nu+1}\frac{1}{jr}\leq \frac{C}{r^{2\nu+2}}\sum_{j=1}^{k [r]} j^{2\nu} \leq \frac{C}{(k[r])^{2\nu+2}} \sum_{j=1}^{k [r]} j^{2\nu}\\
    &\leq \frac{C}{(k[r])^{2\nu+2}}\left\{\begin{array}{ll}
       (k[r])^{2\nu+1},  & \textrm{if }\nu\geq 0, \\
       \displaystyle 1+\int_1^{k[r]} s^{2\nu} ds =1+\frac{(k[r])^{2\nu+1}-1}{2\nu+1}, &  \textrm{if }\nu<0.
    \end{array}\right.
\end{align*}
Hence, 
\[\lim\limits_{r\rightarrow\infty} \|I_{r;2}^k(x,\cdot)\|_{\mathcal H}=0.\]
Next we are going to see that 
\begin{align}\label{4.2}
   \lim_{r\rightarrow \infty}I_{r;1}^k(x,t)=\int_0^{\pi k}{\mathfrak F}_x(y)dy,\quad \mbox{ in }\mathcal{H},
\end{align}
where, as it was shown, the integral can be understood in the $\mathcal{H}$-Riemann sense.

Assume that $\nu>-\frac12$. Recalling that $s_{\nu,j}=\pi \left(j+\frac{\nu}{2}-\frac14+O\left(\frac1j\right)\right)$, for $j\in \mathbb N$, it follows that
\[\frac{s_{\nu, k[r]}}{r}=\frac{\pi \left(k[r]+\frac{\nu}{2}-\frac14+O\left(\frac{1}{k[r]}\right)\right)}{r}=\pi k +\frac{\pi \left(k([r]-r)+\frac{\nu}{2}-\frac14+O\left(\frac{1}{k[r]}\right)\right)}{r}\underset{r\rightarrow \infty}{\longrightarrow} \pi k.\]

Therefore, there exists $\eta>0$ such that $\displaystyle\frac{s_{\nu, k[r]}}{r} <\pi k+1$, for every $r>\eta$.

We distinguish three cases. Firstly, we suppose that $\displaystyle\frac{s_{\nu, k[r]+1}}{r}<\pi k$. If $j\in \mathbb N$ is such that $\displaystyle\frac{s_{\nu, k[r]+j}}{r} <\pi k$, we have that $j<-\frac{\nu}{2}+\frac12+k$, when $r>\eta_1$, for some $\eta_1>\eta$. For every $r>\eta_1$, we define $j_r\in \mathbb N$ such that $\displaystyle\frac{s_{\nu, k[r]+j_r}}{r} <\pi k$ and $\displaystyle\frac{s_{\nu, k[r]+j_r+1}}{r} \geq \pi k$. We write
\begin{align*}
     I_{r;1}^k(x,t) &=\sum_{j=1}^{k [r]+j_r} \left(1-\frac{s_{\nu,j}^2}{t^2 r^2}\right)_+^{\alpha-1} \frac{s_{\nu,j}^2}{t^2 r^2} h_{\nu}(f)\left(\frac{s_{\nu,j}}{r}\right)\sqrt{\frac xr s_{\nu,j}} J_\nu\left(\frac xr s_{\nu,j}\right)\frac{\mathfrak{s}_{\nu, j+1}-s_{\nu,j}}{r}\\
   &\quad -\sum_{j=k[r]+1}^{k [r]+j_r} \left(1-\frac{s_{\nu,j}^2}{t^2 r^2}\right)_+^{\alpha-1} \frac{s_{\nu,j}^2}{t^2 r^2} h_{\nu}(f)\left(\frac{s_{\nu,j}}{r}\right)\sqrt{\frac xr s_{\nu,j}} J_\nu\left(\frac xr s_{\nu,j}\right)\frac{\mathfrak{s}_{\nu, j+1}-s_{\nu,j}}{r}\\
   &:=A_1(r)(t)+B_1(r)(t), \quad r>\eta_1.
\end{align*}
Here $\mathfrak{s}_{\nu, j+1}=s_{\nu, j+1}$, $j=1,\dots, k[r]+j_r-1$, and $\mathfrak{s}_{\nu, k[r]+j_r+1}=k\pi$. Thus, we have that 
\[\|B_1(r)\|_{\mathcal H}\leq C \sum_{j=k[r]+1}^{k [r]+j_r}\left(\frac{s_{\nu,j}}{r}\right)^{2\nu+1}\frac{\mathfrak{s}_{\nu, j+1}-s_{\nu,j}}{r}\leq \frac Cr, \quad r>\eta_1.\]

Suppose now $\displaystyle\frac{s_{\nu, k[r]+1}}{r} >\pi k$. If $j\in \mathbb N$ and $\displaystyle\frac{s_{\nu, k[r]-j}}{r}>\pi k$, then $j<\max\{1,\nu/2\}$ when $r>\eta_2$ for some $\eta_2>\eta$. We redefine, for this case, $j_r\in \mathbb N$ such that $\displaystyle\frac{s_{\nu, k[r]-j_r}}{r}<\pi k$ and $\displaystyle\frac{s_{\nu, k[r]-j_r+1}}{r}\geq \pi k$, where $r>\eta_2$.

We write 
\begin{align*}
      I_{r;1}^k(x,t)  &=\left[ \sum_{j=1}^{k [r]-j_r-1} \left(1-\frac{s_{\nu,j}^2}{t^2 r^2}\right)_+^{\alpha-1} \frac{s_{\nu,j}^2}{t^2 r^2} h_{\nu}(f)\left(\frac{s_{\nu,j}}{r}\right)\sqrt{\frac xr s_{\nu,j}} J_\nu\left(\frac xr s_{\nu,j}\right)\frac{s_{\nu, j+1}-s_{\nu,j}}{r}\right.\\
    &\quad + \left.\left(1-\frac{s_{\nu,k [r]-j_r}^2}{t^2 r^2}\right)_+^{\alpha-1} \frac{s_{\nu,k [r]-j_r}^2}{t^2 r^2} h_{\nu}(f)\left(\frac{s_{\nu,k [r]-j_r}}{r}\right)\sqrt{\frac xr s_{\nu,k [r]-j_r}} J_\nu\left(\frac xr s_{\nu,k [r]-j_r}\right)\frac{\pi k-s_{\nu,k [r]-j_r}}{r}\right]\\
    &\quad + \left[\left(1-\frac{s_{\nu,k [r]-j_r}^2}{t^2 r^2}\right)_+^{\alpha-1} \frac{s_{\nu,k [r]-j_r}^2}{t^2 r^2} h_{\nu}(f)\left(\frac{s_{\nu,k [r]-j_r}}{r}\right)\sqrt{\frac xr s_{\nu,k [r]-j_r}} J_\nu\left(\frac xr s_{\nu,k [r]-j_r}\right)\frac{s_{\nu,k [r]-j_r+1}-\pi k}{r}\right.\\
    &\quad +\left.\sum_{j=k [r]-j_r+1}^{k [r]} \left(1-\frac{s_{\nu,j}^2}{t^2 r^2}\right)_+^{\alpha-1} \frac{s_{\nu,j}^2}{t^2 r^2} h_{\nu}(f)\left(\frac{s_{\nu,j}}{r}\right)\sqrt{\frac xr s_{\nu,j}} J_\nu\left(\frac xr s_{\nu,j}\right)\frac{s_{\nu, j+1}-s_{\nu,j}}{r}\right]\\
    &=A_2(r)(t)+B_2(r)(t), \quad r>\eta_2.
\end{align*}
We have that 
\[\|B_2(r)\|_{\mathcal H}\leq C \sum_{j=k[r]-j_r+1}^{k [r]}\left(\frac{s_{\nu,j}}{r}\right)^{2\nu+1}\frac{1}{r}\leq \frac Cr, \quad r>\eta_2.\]

For the third case, that is, when $\displaystyle\frac{s_{\nu, k[r]+1}}{r}=\pi k$, we define $A_3(r)(t)=I_{r;1}^k(x,t)$.

We observe that $A_j(r)$, $j=1,2,3$, are Riemann sums for ${\mathfrak F}_x$ on $[0,\pi k]$ with amplitude smaller than $C/r$, where $C>0$ does not depend on $r$. Then, if $\varepsilon >0$, there exists $\eta_3 >0$ such that
$$
\left\|\int_0^{\pi k}\mathfrak{F}_x(y)dy-A_j(r)\right\|_{\mathcal{H}}<\varepsilon,\quad r>\eta_3\mbox{ and }j=1,2,3.
$$
We conclude that if $\varepsilon >0$, there exists $\eta >0$ such that 
$$
\left\|\int_0^{\pi k}\mathfrak{F}_x(y)dy-I_{r;1}^k(x,\cdot)\right\|_{\mathcal{H}}<\varepsilon,\quad r >\eta .
$$

Assume now that $-1<\nu\leq -\frac12$. In this case we must do some more work on this because the function $z^{\frac{1}{2}}J_\nu(z)$ is not bounded in any neighbourhood of the origin.

Let $\varepsilon >0$. According to \eqref{4.3} we can find $\lambda \in (0,\pi k)$ such that 
$$
\int_0^\lambda \left\|\mathfrak{F}_x(y)\right\|_{\mathcal{H}}dy\leq C\int_0^\lambda y^{2\nu+1}dy<\varepsilon.
$$
We choose $r_1>0$ such that $\displaystyle\frac{s_{\nu ,j+1}-s_{\nu ,j}}{r}<\frac{\lambda}{2}$, $j\in \mathbb{N}$ and $r >r _1$.

Let $r >r _1$. There exists $k_r$ such that $\displaystyle\frac{s_{\nu ,k_r}}{r} \leq \frac{\lambda}{2}$ and $\displaystyle\frac{s_{\nu,k_{r+1}}}{r} > \frac{\lambda}{2}$. Also, there exists $M>0$ such that
$$
\left(\frac{s_{\nu ,j}}{s_{\nu ,j+1}}\right)^{2\nu +1}\leq M,\quad j\in \mathbb{N}.
$$
Then 
\begin{align*}
\left\|I_{r;1}^k(x,\cdot)\right\|_{\mathcal{H}}&\leq C\sum_{j=1}^{k_r}\left(\frac{s_{\nu ,j}}{r}\right)^{2\nu+1}\frac{s_{\nu ,j+1}-s_{\nu ,j}}{r}\\
& \leq C\sum_{j=1}^{k_r}\left(\frac{s_{\nu ,j+1}}{r}\right)^{2\nu+1}\frac{s_{\nu ,j+1}-s_{\nu ,j}}{r}\leq C\int_0^\lambda y^{2\nu +1}dy\leq C\varepsilon.
\end{align*}

By using the same arguments as in the previous case, we obtain that there exists $\eta >0$ such that
$$
\left\|\int_0^{\pi k}\mathfrak{F}(y)dy-I_{r;1}^k(x,\cdot)\right\|_{\mathcal{H}}<\varepsilon,\quad r >\eta .
$$
Thus \eqref{4.2} is established.

Let $0<a<b<\infty$ and $0<c<d<\infty$. By proceeding as above we get, for every $x\in [a,b]$,
\begin{align*}
\int_c^d\left|I_{r _{\phi (j)}}^{\psi (k)}(x,t)-\left(\int_0^{\pi \psi(k)}\mathfrak{F}_x(y)dy\right)(t)\right|\frac{dt}{t}&\leq \left(\log\frac{d}{c}\right)^{\frac{1}{2}}\left\|I_{r _{\phi (j)}}^{\psi (k)}(x,\cdot)-\left(\int_0^{\pi \psi(k)}\mathfrak{F}_x(y)dy\right)(\cdot)\right\|_{\mathcal{H}}\\
&\hspace{-4cm} \leq \sum_{i=1}^{\psi (k)[r_{\phi (j)}]}\left\|\left(1-\frac{s_{\nu ,i}^2}{t^2r_{\phi(j)} ^2}\right)^{\alpha -1}_{+}\frac{s_{\nu ,i}^2}{t^2r_{\phi(j)} ^2}\right\|_{\mathcal{H}}\left|c_{\nu ,i}(f_{r _{\phi (j)}})\right|\left|\phi _{\nu ,i}\left(\frac{x}{r_{\phi(j)}}\right)\right| \\
&\hspace{-3cm} +\int_0^{\pi \psi (k)}\left\|\left(1-\frac{z^2}{t^2}\right)_+^{\alpha -1}\frac{z^2}{t^2}\right\|_{\mathcal{H}}|h_\nu (f)(z)||\sqrt{xz}J_\nu (xz)|dz\\
&\hspace{-4cm} \leq C\left(\sum_{i=1}^{\psi (k)[r_{\phi (j)}]}\left|c_{\nu ,i}(f_{r _{\phi (j)}})\right|\left|\phi _{\nu ,i}\left(\frac{x}{r_{\phi(j)}}\right)\right|+\int_0^{\pi \psi (k)}|h_\nu (f)(z)||\sqrt{xz}J_\nu (xz)|dz\right)\leq C.
\end{align*}
Dominated Convergence Theorem allows us to obtain
\begin{align*}
\lim_{j\rightarrow \infty }\int_a^b\int_c^d\left|I_{r _{\phi (j)}}^{\psi (k)}(x,t)-\left(\int_0^{\pi \psi(k)}\mathfrak{F}_x(y)dy\right)(t)\right|\frac{dt}{t}dx&\\
&\hspace{-6cm}=\int_a^b\lim_{j\rightarrow \infty}\int_c^d\left|I_{r _{\phi (j)}}^{\psi (k)}(x,t)-\left(\int_0^{\pi \psi(k)}\mathfrak{F}_x(y)dy\right)(t)\right|\frac{dt}{t}dx\\
&\hspace{-6cm}\leq \left(\log \frac{d}{c}\right)^{\frac{1}{2}}\int_a^b\lim_{j\rightarrow \infty}\left(\int_c^d\left|I_{r _{\phi (j)}}^{\psi (k)}(x,t)-\left(\int_0^{\pi \psi(k)}\mathfrak{F}_x(y)dy\right)(t)\right|^2\frac{dt}{t}\right)^{\frac{1}{2}}dx=0.
\end{align*}
Then,
$$
F^k(x,t)=\left(\int_0^{\pi \psi(k)}\mathfrak{F}_x(y)dy\right)(t),\quad \mbox{ a.e. }(x,t)\in (0,\infty )^2.
$$
Hence,
$$
 F(x,t)=\lim_{k\rightarrow \infty}\left(\int_0^{\pi \psi(k)}\mathfrak{F}_x(y)dy\right)(t),\quad \mbox{ a.e. }(x,t)\in (0,\infty )^2.
$$
Let $h\in L^2((0,\infty ),\mathcal{H})$, $0<a<b<+\infty$ and $0<c<d<+\infty$. We have that
\begin{align*}
\int_a^b\int_c^d\left(\int_0^{\pi \psi(k)}\mathfrak{F}_x(y)dy\right)(t)h(x,t)\frac{dt}{t}dx&=\int_a^b\int_0^{\pi \psi (k)}\int_c^d[\mathfrak{F}_x(y)](t)h(x,t)\frac{dt}{t}dydx\\
&=\int_a^b\int_c^d\int_0^{\pi \psi(k)}[\mathfrak{F}_x(y)](t)dyh(x,t)\frac{dt}{t},\quad k\in \mathbb{N}.
\end{align*}
The last equality is justified because, by proceeding as in \eqref{4.3}, we get 
\begin{align*}
\int_a^b\int_0^{\pi \psi (k)}\int_c^d[\mathfrak{F}_x(y)](t)h(x,t)\frac{dt}{t}dydx&\leq
\int_a^b\int_0^{\pi \psi (k)}\|\mathfrak{F}_x(y)\|_{\mathcal{H}}\|h(x,\cdot )\|_{\mathcal{H}}dydx\\
&\leq C\int_a^b\int_0^{\pi \psi (k)}|h_\nu (f)(y)||\sqrt{xy}J_\nu (xy)|dydx\\
&\leq C\int_a^bx^{\nu +\frac{1}{2}}dx\int_0^{\pi \psi (k)}y^{2\nu +1}dy<\infty.
\end{align*}

We deduce that
$$
\left(\int_0^{\pi \psi (k)}\mathfrak{F}_x(y)dy\right)(t)=\int_0^{\pi \psi (k)}\left(1-\frac{y^2}{t^2}\right)_+^{\alpha -1}\frac{y^2}{t^2}h_\nu (f)(y)\sqrt{xy}J_\nu (xy)dy,\quad \mbox{ a.e. }(x,t)\in (0,\infty)^2.
$$

Hence,
$$
F(x,y)=\int_0^\infty\left(1-\frac{y^2}{t^2}\right)_+^{\alpha -1}\frac{y^2}{t^2}h_\nu (f)(y)\sqrt{xy}J_\nu (xy)dy=\frac{t}{2\alpha}\partial_tR_t^{\alpha ,B_\nu}(f)(x),\quad \mbox{ a.e. }(x,t)\in (0,\infty )^2.
$$

We conclude that 
$$
\|G_\alpha ^{B_\nu}f\|_{L^p (0,\infty )}\leq C\|f\|_{L^p(0,\infty )}.
$$
Thus, the proof is finished.

\section{Proof of Theorem \ref{Theorem1.4}}
We consider the operator
$$
m(L_\nu )f=\sum_{j=1}^\infty m_jc_{\nu ,j}(f)\phi _{\nu ,j},\quad f\in L^2(0,1).
$$
Since $\{m_j\}_{j=1}^\infty \in \ell ^\infty$, $m(L_\nu)$ is bounded on $L^2(0,1)$. 

Let $f\in L^2(0,1)$. We have that
$$
G_1^{L_\nu }(f)(x)=\left(\int_0^\infty \left|\frac{\partial}{\partial t}R_t^{1,L_\nu }(f)(x)\right|^2tdt\right)^{\frac{1}{2}}=\left(\int_0^\infty |\mathcal{H}_t^{1, L_\nu}(f)(x)|^2\frac{dt}{t}\right)^{\frac{1}{2}},\quad x\in (0,1),
$$
where
$$
\mathcal{H}_t^{1, L_\nu}(f)(x)=\sum_{j=1}^\infty \chi _{(0,1]}\left(\frac{s_{\nu ,j}}{t}\right)\left(\frac{s_{\nu ,j}}{t}\right)^2c_{\nu ,j}(f)\phi _{\nu ,j}(x),\quad x\in (0,1),\;t>0.
$$
Then,
\begin{align*}
\mathcal{H}_t^{1, L_\nu}(m(L_\nu)f)(x)&=\sum_{j=1}^\infty \chi _{(0,1]}\left(\frac{s_{\nu ,j}}{t}\right)\left(\frac{s_{\nu ,j}}{t}\right)^2m_jc_{\nu ,j}(f)\phi _{\nu ,j}(x)\\
&=\sum_{j=1}^{j(t)} \left(\frac{s_{\nu ,j}}{t}\right)^2m_jc_{\nu,j}(f)\phi _{\nu ,j}(x),\quad x\in (0,1),\;t>0,
\end{align*}
where $j(t)=\max\{j\in \mathbb{N}: s_{\nu ,j}\leq t\}$.

By using Abel summation formula we get
\begin{align*}
\mathcal{H}_t^{1, L_\nu}(m(L_\nu)f)(x)&=m_{j(t)+1}\sum_{j=1}^{j(t)}
\left(\frac{s_{\nu ,j}}{t}\right)^2c_{\nu ,j}(f)\phi _{\nu ,j}(x)-\sum_{k=1}^{j(t)}(m_{k+1}-m_k)\sum_{j=1}^k\left(\frac{s_{\nu ,j}}{t}\right)^2c_{\nu ,j}(f)\phi _{\nu ,j}(x)\\
&=m_{j(t)+1}\mathcal{H}_t^{1, L_\nu}(f)(x)-\sum_{k=1}^{j(t)}(m_{k+1}-m_k)\sum_{j=1}^k\left(\frac{s_{\nu ,j}}{t}\right)^2c_{\nu ,j}(f)\phi _{\nu ,j}(x), \quad x\in (0,1),\;t>0.
\end{align*}
We can write, for each $x\in (0,1)$ and $t>0$,
\begin{align*}
\left|\sum_{k=1}^{j(t)}(m_{k+1}-m_k)\sum_{j=1}^k\left(\frac{s_{\nu ,j}}{t}\right)^2c_{\nu ,j}(f)\phi _{\nu ,j}(x)\right| &\leq \sum_{k=1}^{j(t)}\frac{|m_{k+1}-m_k|}{t^2(s_{\nu, k+1}-s_{\nu ,k})}\int_{s_{\nu,k}}^{s_{\nu ,k+1}}\left|\sum_{j=1}^ks_{\nu,j}^2c_{\nu ,j}(f)\phi _{\nu ,j}(x)\right|du\\
&\hspace{-5cm}=\sum_{k=1}^{j(t)}\frac{|m_{k+1}-m_k|}{t^2(s_{\nu, k+1}-s_{\nu ,k})}\int_{s_{\nu,k}}^{s_{\nu ,k+1}}u^2|\mathcal{H}_u^{1,L_\nu }(f)(x)|du\\
&\hspace{-5cm}\leq \frac{1}{t^2}\left(\sum_{k=1}^{j(t)}|m_{k+1}-m_k|^2\right)^{\frac{1}{2}}
\left(\sum_{k=1}^{j(t)}\frac{1}{(s_{\nu, k+1}-s_{\nu ,k})^2}\int_{s_{\nu,k}}^{s_{\nu ,k+1}}u^5du\int_{s_{\nu,k}}^{s_{\nu ,k+1}}|\mathcal{H}_u^{1,L_\nu }(f)(x)|^2\frac{du}{u}\right)^{\frac{1}{2}}.
\end{align*}

Since $s_{\nu ,k+1}-s_{\nu ,k}=\pi +O(1/k)$, $k\in \mathbb{N}$, it follows that
\begin{align*}
&\left|\sum_{k=1}^{j(t)}(m_{k+1}-m_k)\sum_{j=1}^k\left(\frac{s_{\nu ,j}}{t}\right)^2c_{\nu ,j}(f)\phi _{\nu ,j}(x)\right|\\
&\hspace{5cm}\leq \frac{C}{t^2}\left(\sum_{k=1}^{j(t)}|m_{k+1}-m_k|^2\right)^{\frac{1}{2}}
\left(\sum_{k=1}^{j(t)}\int_{s_{\nu,k}}^{s_{\nu ,k+1}}|\mathcal{H}_u^{1,L_\nu }(f)(x)|^2\frac{du}{u}\right)^{\frac{1}{2}}\\
&\hspace{5cm}\leq \frac{C}{t^2}\left(\sum_{k=1}^{j(t)}|m_{k+1}-m_k|^2\right)^{\frac{1}{2}}
\left(\int_{s_{\nu,1}}^{s_{\nu ,j(t)+1}}|\mathcal{H}_u^{1,L_\nu }(f)(x)|^2\frac{du}{u}\right)^{\frac{1}{2}}.
\end{align*}
Then,
\begin{align*}
\int_0^\infty \left|\sum_{k=1}^{j(t)}(m_{k+1}-m_k)\sum_{j=1}^k\left(\frac{s_{\nu ,j}}{t}\right)^2c_{\nu ,j}(f)\phi _{\nu ,j}(x)\right|^2\frac{dt}{t}& = \int_{s_{\nu ,1}}^\infty \left|\sum_{k=1}^{j(t)}(m_{k+1}-m_k)\sum_{j=1}^k\left(\frac{s_{\nu ,j}}{t}\right)^2c_{\nu ,j}(f)\phi _{\nu ,j}(x)\right|^2\frac{dt}{t} \\
& \hspace{-3cm}\leq C\int_{s_{\nu ,1}}^\infty \frac{1}{t^5}\sum_{k=1}^{j(t)}|m_{k+1}-m_k|^2dt\int_0^\infty |\mathcal{H}_u^{1,L_\nu }(f)(x)|^2\frac{du}{u}\\
&\hspace{-3cm} \leq C\int_0^\infty |\mathcal{H}_u^{1,L_\nu }(f)(x)|^2\frac{du}{u}\sum_{\ell =1}^\infty \int_{s_{\nu ,\ell}}^{s_{\nu ,\ell +1}}\frac{1}{t^5}\sum_{k=1}^{j(t)}|m_{k+1}-m_k|^2dt\\
&\hspace{-3cm} \leq C\int_0^\infty |\mathcal{H}_u^{1,L_\nu }(f)(x)|^2\frac{du}{u}\sum_{\ell =1}^\infty\frac{s_{\nu ,\ell+1}^4-s_{\nu ,\ell }^4}{s_{\nu ,\ell }^4s_{\nu ,\ell +1}^4}\sum_{k=1}^\ell |m_{k+1}-m_k|^2\\
&\hspace{-3cm} \leq C\int_0^\infty |\mathcal{H}_u^{1,L_\nu }(f)(x)|^2\frac{du}{u}\sum_{\ell =1}^\infty\frac{1}{\ell ^5} \sum_{k=1}^\ell |m_{k+1}-m_k|^2,\quad x\in (0,1).
\end{align*}

We conclude that
\begin{equation}\label{5.1}
G_1^{L_\nu }(m(L_\nu)f)(x)\leq |||m|||G_1^{L_\nu}(f)(x),\quad x\in (0,1).
\end{equation}

By using Theorem \ref{Theorem1.1} and \eqref{5.1} the proof of this theorem can be finished.

%\bibliographystyle{acm}
%\bibliography{FourierBesselSteinSquare}

%\def\cprime{$'$} \def\ocirc#1{\ifmmode\setbox0=\hbox{$#1$}\dimen0=\ht0
  %\advance\dimen0 by1pt\rlap{\hbox to\wd0{\hss\raise\dimen0
  %\hbox{\hskip.2em$\scriptscriptstyle\circ$}\hss}}#1\else {\accent"17 #1}\fi}

\end{document}